\documentclass[<option>]{elsarticle}
\usepackage{amsfonts}
\usepackage{tipa}
\usepackage{mathrsfs}
\usepackage{amssymb}
\usepackage{latexsym,bm}
\usepackage{amsthm}
\usepackage{yhmath}
\usepackage{booktabs}
\usepackage{multirow}
\usepackage[center]{caption}
\usepackage{graphics}
\usepackage{subfig}
\usepackage{color}
\usepackage{algorithmic}
\usepackage{algorithm}
\usepackage{geometry}
\usepackage{textcomp}
\journal{arXiv}
\newtheorem{rem}{Remark}[section]
\begin{document}
\begin{frontmatter}
\title{The use of a time-fractional transport model for performing computational homogenisation of 2D heterogeneous media exhibiting memory effects}

\author[els]{Libo~Feng}
\ead{fenglibo2012@126.com}

\author[els,qv]{Ian Turner\corref{cor1}}
\ead{i.turner@qut.edu.au}\cortext[cor1]{Corresponding author.}

\author[rvh]{Patrick Perr\'{e}}
\ead{patrick.perre@centralesupelec.fr}

\author[els,qv,rvt]{Kevin Burrage}
\ead{kevin.burrage@qut.edu.au}

\address[els]{School of Mathematical Sciences, Queensland University of Technology,  GPO Box 2434, Brisbane, QLD. 4001, Australia}
\address[qv]{Australian Research Council Centre of Excellence for Mathematical and Statistical Frontiers (ACEMS), Queensland University of Technology (QUT), Brisbane, Australia}
\address[rvh]{Laboratory of Chemical Engineering and Materials (LGPM), CentraleSup\'{e}lec, Universit\'{e} Paris-Saclay, France}
\address[rvt]{Visiting Professor, Department of Computer Science, University of Oxford, OXI 3QD, UK}

\begin{abstract}
In this work, a two-dimensional time-fractional subdiffusion model is developed to investigate the underlying transport phenomena evolving in a binary medium comprised of two sub-domains occupied by homogeneous material. We utilise an unstructured mesh control volume method to validate the model against a derived semi-analytical solution for a class of two-layered problems. This generalised transport model is then used to perform computational homogenisation on various two-dimensional heterogenous porous media. A key contribution of our work is to extend the classical homogenisation theory to accommodate the new framework and show that the effective diffusivity tensor can be computed once the cell problems reach steady state at the microscopic scale. We verify the theory for binary media via a series of well-known test problems and then investigate media having inclusions that exhibit a molecular relaxation (memory) effect. Finally, we apply the generalised transport model to estimate the bound water diffusivity tensor on cellular structures obtained from environmental scanning electron microscope (ESEM) images for Spruce wood and Australian hardwood. A highlight of our work is that the computed diffusivity for the heterogeneous media with molecular relaxation is quite different from the classical diffusion cases, being dominated at steady-state by the material with memory effects.
\end{abstract}

\end{frontmatter}
\section{Introduction}\label{sec1}

Originating from Fourier's and Darcy's law, numerous real-world problems involve transport processes in heterogeneous media, especially multi-scale transport processes, in which the medium properties differ spatially. One important multi-scale problem is the distributed microstructure model or double/dual-porosity model \cite{Carr14,Lewandowska,Szymkiewicz}, where a large scale domain has a number of small isolated or disconnected inclusions embedded within it.

To model this type of problem, macroscopic and microscopic equations are used to describe the global transport in the connected domain and the local transport in the isolated domains, respectively. Some common approaches for the dual-scale problems include the heterogeneous multi-scale method \cite{Carr14}, distributed microstructure model \cite{Lewandowska} and ``equation-free'' approach \cite{Samaey}. In \cite{Carr16}, Carr et al. proposed an extended distributed microstructure model, which can effectively describe the non-equilibrium field in the inclusions for the water flow in unsaturated soils. However, a number of recent studies highlight that abnormal transport phenomena is evident in some complex heterogeneous media \cite{Turner11,Bueno,perre2019a}, which deviates from the classical Fickian diffusion. Generally, asymptotic long-time behaviour of the mean square displacement is utilised to characterise anomalous diffusion, which has the form
\begin{align*}
\langle x^2(t)\rangle \sim \frac{2K_\gamma}{\Gamma(1+\gamma)}t^\gamma,\quad t\to \infty,
\end{align*}
where $K_\gamma$ is the generalised diffusion coefficient and $\gamma ~(0<\gamma<1)$ is the anomalous diffusion exponent. It turns out that mathematical models involving fractional derivatives are useful tools to treat the anomalous diffusive transport in heterogeneous media \cite{Turner11,Fomin,Ralf2000}, which motivates us to explore these generalised transport equations in this research.

In the seminal paper by Whitaker on the method of volume averaging, macroscopic equations describing the drying process were derived with nonlinear effective parameters within their definitions. This theory also included a set of closure equations that relates these parameters to the pore morphology of the medium, which facilitates their prediction. Homogenization theory can also be used to predict these parameters from the microstructure of the porous medium \cite{Allaire,Hornung}.  The process starts with a full description of the model at the scale of the heterogeneities and upscales the relevant information to derive the macroscopic conservation law by replacing the actual geometry with a periodic idealisation. The macroscopic equation for this homogenised medium is then derived in the limit as the heterogeneities tend to zero. In the classical diffusion model, the effective parameters are defined in terms of the solution of an elliptic equation on the period, subject to periodic boundary conditions (see for example, \cite{Carr14,Lewandowska,Szymkiewicz,Carr16, Bensoussan}). A key contribution of our work is the derivation of a time-evolutionary homogenisation theory that is applied to predict the diffusivity of porous media that exhibit memory effects. A good example of such a material is wood, in which the cellular structure undergoes a molecular relaxation phenomenon during drying \cite{Olek,Olek2}. To model this phenomenon we propose a generalised transport equation involving a time-fractional operator. A challenge is then the way to deal with the fluxes at the interfaces between two different media when one exhibits memory and the other does not. We show how to correctly model this scenario, paying careful attention to the treatment of the interfacial boundary between the different media.

The types of fractional operators utilised in these models  mainly are time-fractional derivatives, which is generally used to describe the algebraic decaying or power-law waiting time in a L\'evy process. For the time-fractional derivative, it has a wide application in many problems \cite{Yuste2001,Benson2003,Langlands2006,Henry2008,Sun2018,Feng2017,Feng2018,Feng2019}. We start with the following two-dimensional time-fractional subdiffusion equation with variable coefficients
\begin{align}\label{eq1}
\frac{\partial u(x,y,t)}{\partial t}={_0^RD^{1-\gamma}_{t}}\left[\nabla\cdot(\bm{Q}\nabla u(x,y,t))\right],\quad(x,y,t)\in \Omega\times (0,T],
\end{align}
subject to the initial and Neumann boundary conditions given, respectively, as
\begin{align} \label{eq2}
u(x,y,0)&=\phi_1(x,y), ~ (x,y)\in\overline{\Omega},\quad {_0^RD^{1-\gamma}_{t}}\left[{\bm{Q}}\nabla u(x,y,t)\cdot {\bm{n}}\right]={\bm{\psi}}(x,y,t), ~ (x,y,t)\in \partial \Omega\times (0,T],
\end{align}
where $\nabla=[\frac{\partial}{\partial x},\frac{\partial}{\partial y}]^T$, $\bm{Q}={\rm{diag}}[q_1(x,y),q_2(x,y)]$, ${\bm{\psi}}(x,y,t)={\rm{diag}}[\psi_1(x,y),\psi_2(x,y)]$, $\bm{n}$ is the unit vector normal to $\partial \Omega$ outward to $\Omega$ and ${_0D^{1-\gamma}_{t}}$ $(0<\gamma<1)$ is the Riemann-Liouville fractional derivative defined by
\begin{align*}
{_0^{R}D_t^{1-\gamma}}u(x,y,t)=\frac{1}{\Gamma(\gamma)}\frac{\partial}{\partial t}\int_0^t\frac{u(x,y,\xi)}{(t-\xi)^{1-\gamma}}d\xi.
\end{align*}
As an application, we will extend the equation to a time-fractional model for simulating transport in heterogeneous systems, which is motivated by the models in \cite{Carr16} and \cite{Zeng19}. In \cite{Carr16}, Carr et al. considered an extended distributed microstructure model for gradient-driven transport on a binary medium without memory. In \cite{Zeng19}, Zeng et al. investigated a coupled system of Caputo time-fractional partial diffusion equations. Different from the models in \cite{Carr16} and \cite{Zeng19}, this paper deals with a time-fractional system based on the Riemann-Liouville fractional derivative on domains including irregularly shaped inclusions with memory. We derive the semi-analytical solution for a class of two-layered problems to validate the proposed numerical method and present an homogenisation theory and numerical calculation, for two-dimensional periodic structures with memory effects. The main contributions of this paper are summarised as follows:
\begin{itemize}
\item[$\bullet$] A two-dimensional time-fractional subdiffusion equation with variable coefficients is considered, in which an unstructured mesh control volume method (CVM) is applied to solve the problem. Given the non-smoothness of the solution, the modified weighted shifted Gr\"unwald-Letnikov (WSGL) formula with starting weights is utilised to deal with the non-smooth problem, which is further extended to solve the time-fractional transport models in heterogeneous media.

\item[$\bullet$] A novel time-fractional transport model on a binary medium (a porous medium comprised of two distinct phases) with regular or irregular-shaped inclusions is simulated using the CVM. To verify the accuracy of the computational model we derive a semi-analytical solution for a class of two-layered problems with quasi-periodic boundary conditions, in which the finite Fourier and Laplace transforms together with a numerical inverse Laplace transform technique are used. In addition, the mass balance equation for the layered medium is also developed.

\item[$\bullet$] Homogenisation theory is extended from a classical diffusion equation to a time-fractional diffusion equation for carrying out computational homogenisation on some well-known test problems and then on media having inclusions that exhibit a molecular relaxation effect. The generalised transport model is also applied to estimate the bound water diffusivity tensor on cellular structures obtained from ESEM images for Spruce wood and Australian hardwood (Eucalyptus pilularis). An important finding is that, unlike the classical case, the fractional-order indices have a significant effect on the mass transfer and the equivalent diffusivity is dominated by the material with memory effects.
\end{itemize}
The structure of this paper is as follows. In Section \ref{sec2}, the modified WSGL formulae with starting weights for the time-fractional operator is proposed and the unstructured mesh control volume method for the two-dimensional subdiffusion problem (\ref{eq1}) is introduced. We show that the proposed computational scheme offers second order convergence in time and space. In Section \ref{sec3}, a time-fractional transport model on a binary medium is proposed and the semi-analytical solution for a class of two-layered problems is derived. In Section \ref{sec4}, homogenisation theory for the time-fractional diffusion equation is developed. A series of numerical examples are considered in Sections \ref{sec5} and \ref{sec6} to verify the homogenisation theory, in which different heterogeneous media morphologies and wood cellular periodic structures are simulated. The numerical results show that the fractional-order indices have a significant effect on the mass transfer for the time-fractional transport equation, which is very different to the classical case.  Finally, some conclusions are drawn in Section \ref{sec7}.

\section{The control volume method  for the subdiffusion equation}\label{sec2}
\subsection{Discretisation for the time-fractional operator}\label{sec2.1}
Firstly, we present the grid partition used to discretise equation (\ref{eq1}) in the temporal direction. Define $t_n=n\tau$, $n=0,1,2,\ldots,N$, where $\tau=\frac{T}{N}$ is the uniform temporal step. The Riemann-Liouville time-fractional derivative can be rewritten in a more general form as ${_0^{R}D_t^{\alpha}}v(t)=\frac{1}{\Gamma(1-\alpha)}\frac{d}{d t}\int_0^t\frac{v(\xi)}{(t-\xi)^\alpha}d\xi$, $0<\alpha<1$.
To approximate the Riemann-Liouville fractional derivative, the shifted Gr\"unwald-Letnikov difference operator was proposed \cite{MT04,Wang14}:
$\mathcal{A}_{\tau,p}^{\alpha}v(t)=\tau^{-\alpha}\sum_{k=0}^{\infty}g_k^{(\alpha)}v(t-(k-p)\tau)$,
where $g_k^{(\alpha)}=(-1)^k\binom{\alpha}{k}$ and $p$ is an integer. To improve the accuracy of the approximation formula, some high order schemes based on the WSGL formula have been developed \cite{Wang14,Tian15}. Here we focus on the second order WSGL formula, which has the form
\begin{align*}
\mathcal{B}_{\tau,p,q}^{\alpha}v(t)=\frac{\alpha-2q}{2(p-q)}\mathcal{A}_{\tau,p}^{\alpha}v(t)
+\frac{2p-\alpha}{2(p-q)}\mathcal{A}_{\tau,q}^{\alpha}v(t),
\end{align*}
where $p$ and $q$ are two non equal integers. According to \cite{Tian15}, $p$ and $q$ should be chosen satisfying $|p|\leq 1$ and $|q|\leq 1$. When $(p,q)=(0,-1)$, we have $\frac{\alpha-2q}{2(p-q)}=\frac{2+\alpha}{2}$, $\frac{2p-\alpha}{2(p-q)}=-\frac{\alpha}{2}$. Then at $t=t_n$, we obtain
\begin{align*}
{_{0}^{R}D_t^{\alpha}}v(t_n)=\mathcal{B}_{\tau,p,q}^{\alpha,n}v+O(\tau^2)
=\tau^{-\alpha}\sum_{k=0}^{n}\omega_{n-k}^{(\alpha)}v(t_k)+O(\tau^2),
\end{align*}
where $\mathcal{B}_{\tau,p,q}^{\alpha,n}v=\frac{2+\alpha}{2}\mathcal{A}_{\tau,p}^{\alpha}v(t_n)-
\frac{\alpha}{2}\mathcal{A}_{\tau,q}^{\alpha}v(t_n)$ and the convolution weights $\omega_{k}^{(\alpha)}$ are defined as
$\omega_{0}^{(\alpha)}=\frac{2+\alpha}{2}g_{0}^{(\alpha)}$, $\omega_{k}^{(\alpha)}=\frac{2+\alpha}{2}g_{k}^{(\alpha)}-\frac{\alpha}{2}g_{k-1}^{(\alpha)}$, $k\geq 1$.
Define $\mathcal{D}_{\tau,p,q}^{\alpha,n}v=\frac{1}{2}(\mathcal{B}_{\tau,p,q}^{\alpha,n}v+\mathcal{B}_{\tau,p,q}^{\alpha,n-1}v)$, then at $t=t_{n-\frac{1}{2}}$, we have
\begin{align} \label{eq3}
{_{0}^{R}D_t^{\alpha}}v(t_{n-\frac{1}{2}})=\mathcal{D}_{\tau,p,q}^{\alpha,n}v+O(\tau^2)=
\tau^{-\alpha}\sum_{k=0}^{n}D_{n-k}^{(\alpha)}v(t_k)+O(\tau^2),
\end{align}
where $D_{0}^{(\alpha)}=\frac{\omega_{0}^{(\alpha)}}{2}$, $D_{k}^{(\alpha)}=\frac{\omega_{k}^{(\alpha)}+\omega_{k-1}^{(\alpha)}}{2}$, $k\geq 1$.

\begin{rem}
To guarantee the second-order convergence of (\ref{eq3}), the required conditions are $v(t)\in L_1(\mathbb{R})$ and ${_{0}^{R}D_t^{\alpha+2}}v(t)$ and its Fourier transform belongs to $L_1(\mathbb{R})$ \cite{Tian15}.
\end{rem}
\subsection{The implementation of the control volume method}\label{sec2.2}

\begin{figure}[H]
\begin{center}
\scalebox{0.4}[0.4]{\includegraphics{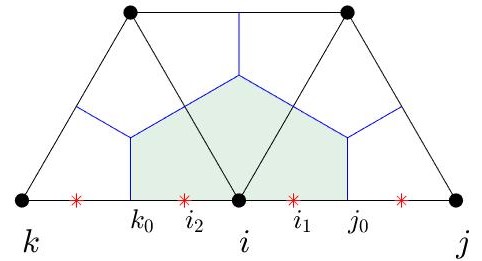}}
\caption{An illustration of a point with control volume on the boundary.}
\label{Fig1}
\end{center}
\end{figure}

To implement the control volume method, we partition the solution domain $\Omega$ into a mesh comprised of non-overlapping triangles. Denote $\mathcal{T}_h$ as the triangulation, $N_e$ the number of triangles and $h$ the maximum diameter of the triangular elements. For a group of triangles with the same vertex, a control volume can be formed by joining the midpoints and barycenters of each triangle. Integrating (\ref{eq1}) over each control volume $V_i$ $(i=1,2,\ldots,N_p)$ and applying Gauss's Divergence Theorem, yields
\begin{align*}
\int_{V_i}\frac{\partial u}{\partial t}dV_i={_0^RD^{1-\gamma}_{t}}\oint_{\partial V_i}(\bm{Q}\nabla u)\cdot{\mathbf{n}}\,d\Gamma_i,
\end{align*}
where $d\Gamma_i$ is an infinitesimal small line segment on $\partial V_i$. Utilising a lumped mass approach to approximate the time derivative, gives
\begin{align} \label{eq4}
\triangle{V_i}\frac{\partial u_i}{\partial t}={_0^RD^{1-\gamma}_{t}}\oint_{\partial V_i}(\bm{Q}\nabla u)\cdot{\mathbf{n}}\,d\Gamma_i,\quad i=1,\ldots,N_p,
\end{align}
where $\triangle{V_i}$ is the volume of the control volume. We assume that the integration path is anticlockwise and the triangle vertices are numbered counter-clockwise. We choose piecewise linear polynomials on the triangles. Within an element triangle $e_p$, $p=1,2,\ldots,N_e$, the shape function can be defined in the form $N_i(x,y)=\frac{1}{2\Delta_{e_p}}\left(a_ix+b_iy+c_i\right)$, $i=1,2,3$, where $a_i$, $b_i$ and $c_i$ are some constants and $\Delta_{e_p}$ is the area of triangle element $p$. Combining these shape functions, we can construct the basis function $l_k(x,y)$ and approximate $u(x,y,t)$ as $u(x,y,t)\approx u_h(x,y,t)=\sum_{k=1}^{N_P}u_k\,l_k(x,y)$. The line integral in (\ref{eq4}) can be approximated by the midpoint formula at each control surface:
\begin{align}
&\oint_{\partial V_i}(\bm{Q}\nabla u)\cdot{\mathbf{n}}\,d\Gamma_i=\oint_{\Gamma_i}q_1(x,y,t)\frac{\partial u}{\partial x}dy-\oint_{\Gamma_i}q_2(x,y,t)\frac{\partial u}{\partial y}dx\nonumber\\
=&\sum_{j=1}^{m_i}\sum_{r=1}^{2}\left( q_1(x,y,t)\frac{\partial u}{\partial x}\right)\bigg|_{(x_r,y_r)}\Delta y_{j,r}^i
-\sum_{j=1}^{m_i}\sum_{r=1}^{2}\left(q_2(x,y,t)\frac{\partial u}{\partial y} \right)\bigg|_{(x_r,y_r)}\Delta x_{j,r}^i\nonumber\\ \label{eq5}
=&\sum_{k=1}^{N_P}\sum_{j=1}^{m_i}\sum_{r=1}^{2}u_k\left( q_1(x,y,t)\frac{\partial l_k(x,y)}{\partial x}\right)\bigg|_{(x_r,y_r)}\Delta y_{j,r}^i
-\sum_{k=1}^{N_P}\sum_{j=1}^{m_i}\sum_{r=1}^{2}u_k\left(q_2(x,y,t)\frac{\partial l_k(x,y)}{\partial y} \right)\bigg|_{(x_r,y_r)}\Delta x_{j,r}^i,
\end{align}
where $(x_r,y_r)$ is the mid-point of the control face and $m_i$ is the number of sub-control volumes associated with the node $i$. For more details, the readers can refer to \cite{Feng19}. Then from (\ref{eq4}), we can derive the following ODE system
\begin{align} \label{eq6}
\bm{M}\frac{d\bm{u}}{dt}={_0^RD^{1-\gamma}_{t}}\bm{Ku}+\bm{F}_b,
\end{align}
where $\bm{M}=\rm{diag}[\triangle V_1, \triangle V_2,\ldots,\triangle V_{N_p}]$, $\bm{u}=[u_1,u_2,\ldots,u_{N_p}]^T$, $\bm{K}$ is the stiffness matrix derived from (\ref{eq5}) and $\bm{F}_b$ is the contribution from the boundary conditions. Now we discuss the treatment of the boundary conditions. For example, a point $i$ is on the boundary with a control volume and $k$ and $j$ are its adjacent points on the boundary (see Figure \ref{Fig1}); $k_0$ is the midpoint of $k$ and $i$ and $j_0$ is the midpoint of $j$ and $i$; $i_1$ is the midpoint of $j_0$ and $i$ and $i_2$ is the midpoint of $k_0$ and $i$. Then according to (\ref{eq2}) and (\ref{eq5}), we obtain
\begin{align*}
\bm{F}_b(i)&={_0^RD^{1-\gamma}_{t}}\left[\left( q_1(x,y,t)\frac{\partial u}{\partial x}\right)\bigg|_{(x_{i_1},y_{i_1})}(y_{j_0}-y_i)
-\left(q_2(x,y,t)\frac{\partial u}{\partial y} \right)\bigg|_{(x_{i_1},y_{i_1})}(x_{j_0}-x_i)\right]\\
&+{_0^RD^{1-\gamma}_{t}}\left[\left( q_1(x,y,t)\frac{\partial u}{\partial x}\right)\bigg|_{(x_{i_2},y_{i_2})}(y_i-y_{k_0})
-\left(q_2(x,y,t)\frac{\partial u}{\partial y} \right)\bigg|_{(x_{i_2},y_{i_2})}(x_i-x_{k_0})\right]\\
&=\psi_1(x_{i_1},y_{i_1},t)(y_{j_0}-y_i)-\psi_2(x_{i_1},y_{i_1},t)(x_{j_0}-x_i)
+\psi_1(x_{i_2},y_{i_2},t)(y_i-y_{k_0})-\psi_1(x_{i_2},y_{i_2},t)(x_i-x_{k_0}).
\end{align*}
At $t=t_{n-\frac{1}{2}}$, applying the Crank-Nicolson scheme to the integer time derivative and the WSGL formula to the time-fractional derivative, we obtain
\begin{align} \label{eq7}
&\bm{M}\frac{\bm{u}^n-\bm{u}^{n-1}}{\tau}
=\tau^{\gamma-1}\sum_{k=1}^{n}D_{n-k}^{(1-\gamma)}\bm{K}(\bm{u}^k-\bm{u}^0)
+\bm{F}_b^{n-\frac{1}{2}}+\frac{t_{n-\frac{1}{2}}^{\gamma-1} \bm{K}\bm{u}^0}{\Gamma(\gamma)}.
\end{align}
\vspace{-4mm}
\vspace{-8mm}
\subsection{Modified numerical scheme with correction terms}\label{sec2.4}
Since the WSGL formulae require some restrictions on the solution (${_{0}^{R}D_t^{\alpha+2}}v(t)\in L_1(\mathbb{R})$), they may not exhibit the expected high order accuracy when the solution is not smooth enough, especially if it involves terms of the form $t^{\sigma_r}$. To overcome this issue, the modified WSGL formulae with correction terms were proposed to deal with problems exhibiting a non-smooth solution \cite{Zeng17,Liu2019}. At $t=t_{n-\frac{1}{2}}$, we have
\begin{align} \label{eq8}
{_{0}^{R}D_t^{\alpha}}v(t_{n-\frac{1}{2}})=\tau^{-\alpha}\sum_{k=0}^{n}D_{n-k}^{(\alpha)}v(t_k)+\tau^{-\alpha}\sum_{k=1}^{m}E_k^{(n,\alpha)}v(t_k)+O(\tau^2),
\end{align}
where the starting weights $E_k^{(n,\alpha)}$ are chosen such that  (\ref{eq8}) is exact for $v(t)=t^{\sigma_r}$ ($1\leq r \leq m$), which leads to the following linear system:
\begin{align*}
\sum_{k=1}^{m}E_k^{(n,\alpha)}k^{\sigma_r}=\frac{\Gamma(\sigma_r+1)}{\Gamma(\sigma_r+1-\alpha)}\left(n-\frac{1}{2}\right)^{\sigma_r-\alpha}
-\sum_{k=0}^{n}D_{n-k}^{(\alpha)}k^{\sigma_r},\quad r=1,\ldots,m.
\end{align*}
In addition, we also need to define the starting weights for the first order time derivative at $t=t_{n-\frac{1}{2}}$, i.e.,
\begin{align}\label{eq9}
\frac{dv(t_{n-\frac{1}{2}})}{dt}=\frac{v(t_n)-v(t_{n-1})}{\tau}+\tau^{-1}\sum_{k=1}^{m}P_k^{(n,\alpha)}v(t_k)+O(\tau^2),
\end{align}
where the starting weights $P_k^{(n,\alpha)}$ are chosen such that (\ref{eq9}) is exact for $v(t)=t^{\sigma_r}$ ($1\leq r \leq m$), which can be obtained from the following equations:
\begin{align*}
\sum_{k=1}^{m}P_k^{(n,\alpha)}k^{\sigma_r}=\sigma_r\left(n-\frac{1}{2}\right)^{\sigma_r-1}-\big(n^{\sigma_r}-(n-1)^{\sigma_r}\big),\quad r=1,\ldots,m.
\end{align*}
In this paper, we choose the correction terms $\sigma_r=r\alpha$, $r=1,2,\ldots,m$, such that $(m+1)\alpha\geq2$. For more details, the readers can refer to \cite{Zeng17}. Then the modified numerical scheme of (\ref{eq7}) with correction terms is
\begin{align}
&\bm{M}\frac{\bm{u}^n-\bm{u}^{n-1}}{\tau}+\frac{1}{\tau}\sum_{k=1}^{m}P_k^{(n,\alpha)}\bm{M}(\bm{u}^k-\bm{u}^0)\nonumber\\\label{eq10}
=&\tau^{\gamma-1}\left[\sum_{k=1}^{n}D_{n-k}^{(1-\gamma)}\bm{K}(\bm{u}^k-\bm{u}^0)
+\sum_{k=1}^{m}E_k^{(n,\alpha)}\bm{K}(\bm{u}^k-\bm{u}^0)\right]+\bm{F}_b^{n-\frac{1}{2}}+\frac{t_{n-\frac{1}{2}}^{\gamma-1} \bm{K}\bm{u}^0}{\Gamma(\gamma)},
\end{align}
which also can be extended to solve the time-fractional transport model (\ref{eq11}).

\subsection{Investigation of the accuracy and convergence order of CVM with the modified WSGL formulae}
We complete this section by assessing the accuracy and convergence order of the CVM method with the modified WSGL formulae for the two-dimensional time-fractional subdiffusion model (\ref{eq1}) on a square domain $\Omega=[0,1] \times [0,1]$. We treat the problem as isotropic $\mathbf{Q}= \tfrac{1}{2} \, \mathbf{I}$ with initial condition $u(x,y,0)=\sin x\sin y$ and boundary conditions ${_0^RD^{1-\gamma}_{t}}\left[\mathbf{Q}\frac{\partial v(x,y,t)}{\partial x} \cdot \mathbf{n} \right]
=t^{\gamma-1}E_{\gamma,\gamma}(-t^{\gamma})\cos x\sin y,(x,y,t)\in \partial \Omega\times (0,T]$, where the Mittag-Leffler function is $E_{\alpha,\beta}(z)=\sum_{k=0}^{\infty} \frac{z^k}{\Gamma(k\alpha+\beta)}$. The exact solution of this problem is $u(x,y,t)=E_{\gamma,1}(-t^{\gamma})\sin x\sin y$.

\begin{table}[H]
\begin{center}
\caption{The error and convergence order of numerical schemes (\ref{eq7}) and (\ref{eq10}) for problem (\ref{eq1}) for different $h$ and $\gamma_1=\gamma_2=\gamma$ with $\tau=1\times10^{-3}$ at $t=1$, in which the number of correction terms $m$ for $\gamma=0.5$ and $\gamma=0.8$ are $m=3$ and $m=2$, respectively.}
\label{tab1}
\begin{tabular}{ccccc}
\toprule
& \multicolumn{2}{c}{No correction} & \multicolumn{2}{c}{Apply correction}\\
\midrule
 $h~(\gamma=0.5)$   &  error & Order & error & Order \\
\midrule
2.6170E-01  & 1.2688E-03  &  --   & 1.3558E-03 &  --   \\
1.5176E-01  & 1.9206E-03  & -0.76 & 4.4375E-04 & 2.05 \\
8.4435E-02  & 2.1995E-03  & -0.23 & 1.0488E-04 & 2.46 \\
4.1624E-02  & 2.2646E-03  & -0.04 & 2.8136E-05 & 1.86  \\
2.1733E-02  & 2.2821E-03  & -0.01 & 7.3674E-06 & 2.06 \\
\midrule
 $h~(\gamma=0.8)$   &  error & Order & error & Order \\
\midrule
2.6170E-01  & 1.2909E-03  &  --   & 1.3607E-03 &  --   \\
1.5176E-01  & 3.7649E-04  & 2.26  & 4.4321E-04 & 2.06 \\
8.4435E-02  & 5.3879E-05  & 3.32  & 1.0482E-04 & 2.46 \\
4.1624E-02  & 5.6789E-05  & -0.07 & 2.8160E-05 & 1.86  \\
2.1733E-02  & 7.2710E-05  & -0.38 & 7.4373E-06 & 2.05 \\
\bottomrule
\end{tabular}
\end{center}
\end{table}
\vspace{-5mm}
We apply the numerical scheme (\ref{eq7}) and the modified numerical scheme (\ref{eq10}) to solve (\ref{eq1}) and present the error and convergence order for different $h$ and $\gamma$ with $\tau=1\times10^{-3}$ at $t=1$ in Table \ref{tab1}. It is straightforward to observe that the regularity of the exact solution is low particularly when $\gamma$ is small, which can cause difficulty for the numerical method to attain the theoretical temporal convergence order. This finding is evident in the results presented in Table \ref{tab1}, where we can see that the numerical scheme (\ref{eq7}) without correction terms fails to provide good convergence behaviour. However, when correction terms are added, the singularity of the solution can be captured effectively and global second-order convergence is achieved. We conclude that the modified numerical scheme (\ref{eq10}) with correction terms is a very effective computational method, thus justifying its choice for solving the generalised transport models considered throughout the following sections.

\section{Application to a heterogeneous medium}\label{sec3}
\subsection{A time-fractional transport model}\label{sec3.1}
\begin{figure}[H]
\begin{center}
\scalebox{0.35}[0.35]{\includegraphics{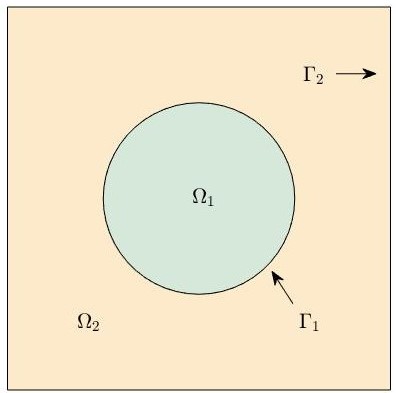}}
\caption{An illustration of a binary medium comprised of two sub-domains occupied by homogeneous material. The square domain (medium 2) is connected and the circular inclusion (medium 1) is isolated or disconnected.}
\label{Fig2}
\end{center}
\end{figure}

In this section, we extend the control volume method to solve a time-fractional transport model on a binary medium consisting of a connected phase $\Omega_2$ and an inclusion (disconnected phase) $\Omega_1$. The time-fractional indices are different for the two media (see Figure \ref{Fig2}). Here we consider the following system:
\begin{equation}\label{eq11}
\left\{\begin{array}{l}
\frac{\partial u(x,y,t)}{\partial t}={_0^RD^{1-\gamma_1}_{t}}\left[\nabla\cdot(\bm{Q}_1\nabla u(x,y,t))\right],\quad(x,y,t)\in \Omega_1\times (0,T],\\
\frac{\partial v(x,y,t)}{\partial t}={_0^RD^{1-\gamma_2}_{t}}\left[\nabla\cdot(\bm{Q}_2\nabla v(x,y,t))\right],\quad(x,y,t)\in \Omega_2\times (0,T],
\end{array}\right.
\end{equation}
\noindent subject to
\begin{align*}
&u(x,y,0)=\phi_1(x,y), \quad (x,y)\in \Omega_1\cup\Gamma_1,\quad
v(x,y,0)=\phi_2(x,y), \quad (x,y)\in \Omega_2\cup\Gamma_1\cup\Gamma_2,\\
&{_0^RD^{1-\gamma_2}_{t}}\left[{\bm{Q}_2}\nabla v(x,y,t)\cdot {\bm{n}}\right]={\bm{\psi}(x,y,t)}, \quad (x,y,t)\in \Gamma_2\times (0,T].
\end{align*}
In addition, we also need the following boundary conditions at the interface $\Gamma_1$ between the two media
\begin{align}\label{eq12}
 u(x,y,t)= v(x,y,t), \quad {_0^RD^{1-\gamma_1}_{t}}(\bm{Q}_1\nabla u(x,y,t)\cdot {\bm{n}})={_0^RD^{1-\gamma_2}_{t}}(\bm{Q}_2\nabla v(x,y,t)\cdot {\bm{n}}),\quad (x,y,t)\in \Gamma_1\times (0,T].
\end{align}
Integrating (\ref{eq11}) over $V_i\in \Omega_1$ and $V_j\in \Omega_2$, respectively, yields
\begin{equation*}
\left\{\begin{array}{l}
\int_{V_i}\frac{\partial u}{\partial t}dV_i={_0^RD^{1-\gamma_1}_{t}}\oint_{\partial V_i}(\bm{Q}_1\nabla u)\cdot{\mathbf{n}}\,d\Gamma_i,\\
\int_{V_j}\frac{\partial v}{\partial t}dV_j={_0^RD^{1-\gamma_2}_{t}}\oint_{\partial V_j}(\bm{Q}_2\nabla v)\cdot{\mathbf{n}}\,d\Gamma_j.
\end{array}\right.
\end{equation*}
\vspace{-6mm}
\noindent Furthermore, the matrix form can be derived as
\begin{equation}\label{eq13}
\left\{\begin{array}{l}
\bm{M_1}\frac{d\bm{u}}{dt}={_0^RD^{1-\gamma_1}_{t}}\bm{K_1u}+ \bm{F}_{b1},\\
\bm{M_2}\frac{d\bm{v}}{dt}={_0^RD^{1-\gamma_2}_{t}}\bm{K_2v}+\bm{F}_{b2}-\bm{F}_{b1},
\end{array}\right.
\end{equation}
where $\bm{u}=[\bm{u_1},\bm{u_{12}}]^T$, $\bm{v}=[\bm{v_1},\bm{v_{12}}]^T$, $\bm{u_1}\in \Omega_1$, $\bm{u_{12}}=\bm{v_{12}}\in\Gamma_1$ (which is the shared part at the interface $\Gamma_1$) and $\bm{v_1}\in \Omega_2$, $\bm{F}_{b1}$ is the contribution from the boundary conditions at the internal interface $\Gamma_1$ and $\bm{F}_{b2}$ is the contribution from the boundary conditions at the external interface $\Gamma_2$. Denote $\bm{U}=[\bm{u_1},\bm{u_{12}},\bm{v_1}]^T$, then the system (\ref{eq13}) can be recast as
\begin{align*}
\bm{M}\frac{d\bm{U}}{dt}={_0^RD^{1-\gamma}_{t}} \bm{KU}+\bm{F}_{b2},
\end{align*}
which can be solved using the same technique discussed in the previous section.
\subsection{The semi-analytical solution for a class of two-layered problems}\label{sec3.2}

\begin{figure}[H]
\begin{center}
\scalebox{0.55}[0.55]{\includegraphics{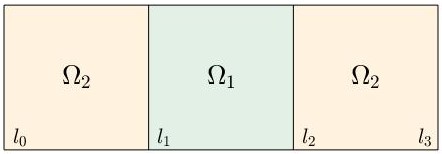}}
\caption{An illustration of a two-layered medium.}
\label{Fig3}
\end{center}
\end{figure}

In this section, we consider the semi-analytical solution for the following 2D time-fractional two-layered problem (see Figure \ref{Fig3}).
\begin{equation}\label{eq14}
\left\{\begin{array}{l}
\frac{\partial u(x,y,t)}{\partial t}={_0^RD^{1-\alpha_1}_{t}}\left[D_{1}\Delta u(x,y,t)\right],\quad(x,y,t)\in \Omega_1\times (0,T],\\
\frac{\partial v(x,y,t)}{\partial t}={_0^RD^{1-\alpha_2}_{t}}\left[D_{2}\Delta  v(x,y,t)\right],\quad(x,y,t)\in \Omega_2\times (0,T],
\end{array}\right.
\end{equation}
where $\Omega_1=(l_1,l_2)\times (0,l_y)$ and $\Omega_2=(l_0,l_3)\times (0,l_y)\backslash \Omega_1$. Since the medium is homogeneous in the $y$-direction and we impose no flux boundary conditions at $y=0$ and $y=l_y$, solving problem (\ref{eq14}) can be reduced to solving the following 1D two-layered problem:
\begin{equation}\label{eq15}
\left\{\begin{array}{l}
\frac{\partial u(x,t)}{\partial t}={_0^RD^{1-\alpha_1}_{t}}\left(D_{1}\frac{\partial^2u(x,t)}{\partial x^2} \right),\quad(x,t)\in (l_1,l_2)\times (0,T],\\
\frac{\partial v(x,t)}{\partial t}={_0^RD^{1-\alpha_2}_{t}}\left(D_{2} \frac{\partial^2v(x,t)}{\partial x^2}\right),\quad(x,t)\in (l_0,l_1)\cup(l_2,l_3)\times (0,T].
\end{array}\right.
\end{equation}
\vspace{-3mm}
\vspace{-6mm}
\noindent In the following, we will consider the equivalent form of (\ref{eq15}), which can be written as
\begin{align} \label{eq16}
\frac{\partial X_i(x,t)}{\partial t}={_0^RD^{1-\gamma_i}_{t}}\left(D_{bi}\frac{\partial^2X_i(x,t)}{\partial x^2} \right),\quad x\in(l_{i-1},l_i), \quad i=1,2,3,
\end{align}
with the initial conditions $X_i(x,0)=X_{i0}(x)$ and quasi-periodic boundary conditions of the following form
\begin{align*}
&{_0^RD^{1-\gamma_1}_{t}}{\left(D_{b1}\frac{\partial X_1}{\partial x}(l_0,t) \right)}={_0^RD^{1-\gamma_3}_{t}}{\left(D_{b3}\frac{\partial X_3}{\partial x}(l_3,t) \right)},\quad X_3(l_3,t)=X_1(l_0,t)+q_0,\nonumber\\
&{_0^RD^{1-\gamma_i}_{t}}{\left(D_{b,i}\frac{\partial X_{i}}{\partial x}(l_i,t) \right)}=
{_0^RD^{1-\gamma_{i+1}}_{t}}{\left(D_{b,i+1}\frac{\partial X_{i+1}}{\partial x}(l_i,t) \right)},\quad X_{i}(l_i,t)=X_{i+1}(l_i,t),\quad i=1,2,
\end{align*}
where $u=X_2$, $v=X_1\cup X_3$, $q_0$ is a constant, $D_{b2}=D_1$, $D_{b1}=D_{b3}=D_2$, $\gamma_2=\alpha_1$ and $\gamma_1=\gamma_3=\alpha_2$. These boundary conditions have been chosen to validate the time evolutionary homogenisation theory proposed in the proceeding section. The problem (\ref{eq16}) will be solved using a combination of finite Fourier and Laplace transforms, which is an extension of the semi-analytical approach discussed in \cite{Carr16b}. First define $d_i=l_i-l_{i-1},i=1,2,3$ and then apply the finite Fourier transform in each layer by noting that the transform associated with the following Sturm-Liouville system in $\varphi_i, i=1,2,3$:
\begin{align*}
-\frac{d^2}{dx^2} \varphi_i = \lambda_i^2 \varphi_i, \quad
\frac{d \varphi_i}{dx}(l_{i-1})=0 ~\mathrm{and}~\frac{d \varphi_i}{dx}(l_i)=0,
\end{align*}
is given by $\widetilde{X}_i(\lambda_{i,m},t) := \langle X_i,\varphi_{i,m} \rangle= \int_{l_{i-1}}^{l_i} \, X_i(x,t)\, \varphi_{i,m}(x)\, dx$, where the eigenvalues $\lambda_{i,m}^2=\frac{m^2\pi^2}{d_i^2}, m\geq0$. The corresponding eigenfunctions are $\varphi_{i,m}(x)=\frac{1}{\sqrt{d_i}}$ for $m=0$ and $\varphi_{i,m}(x)=\sqrt{\frac{2}{{d_i}}}\cos {\left[ \lambda_m (x-l_{i-1})\right]}$ for $m>0$. We introduce the time fractional potential for each layer as $\mathcal{X}_i={_{0}^{R}D}_{t}^{1-\gamma_i}X_i,~i=1,2,3$ and then apply the finite Fourier transform to each system within the $i^{\mathrm{th}}$ layer to obtain
\begin{align*}
\left\langle \frac{\partial X_i}{\partial t}, \varphi_{i,m} \right\rangle &= D_{bi}  \left\langle  \frac{\partial^2 \mathcal{X}_i}{\partial x^2}, \varphi_{i,m} \right\rangle.
\end{align*}
The term on the righthand side is integrated by parts to obtain
\begin{align}
D_{bi} \left\langle \frac{\partial^2 \mathcal{X}_i}{\partial x^2}, \varphi_{i,m} \right\rangle = -D_{bi} \left\{ \left\langle -\frac{d^2}{dx^2} \varphi_{i,m}, \mathcal{X}_i \right\rangle + \left[-  \frac{\partial \mathcal{X}_i}{\partial x}(l_i,t) \varphi_{i,m}(l_i) + \frac{\partial \mathcal{X}_i}{\partial x}(l_{i-1},t) \varphi_{i,m}(l_{i-1})\right] \right\}. \label{eq17}
\end{align}
Next, define the interfacial flux terms as $D_{b1}\frac{\partial \mathcal{X}_1}{\partial x}(l_1,t) = D_{b2} \frac{\partial \mathcal{X}_2}{\partial x}(l_1,t)= v_{12}(t)$, $D_{b1}\frac{\partial \mathcal{X}_1}{\partial x}(l_0,t) = D_{b3} \frac{\partial \mathcal{X}_3}{\partial x}(l_3,t)= v_{13}(t)$, $D_{b2}\frac{\partial \mathcal{X}_2}{\partial x}(l_2,t) = D_{b3} \frac{\partial \mathcal{X}_3}{\partial x}(l_2,t)= v_{23}(t)$. Substituting the relevant boundary condition information into (\ref{eq17}), the transformed layer equations are then given by
\begin{align*}
\frac{\partial \widetilde{X}_1}{\partial t} &=-D_{b1} \lambda_{1,m}^2 \widetilde{\mathcal{X}}_1 + \left[v_{12}(t)\varphi_{1,m}(l_1) -v_{13}(t)\varphi_{1,m}(l_{0}) \right],\\
\frac{\partial \widetilde{X}_2}{\partial t} &=-D_{b2} \lambda_{2,m}^2 \widetilde{\mathcal{X}}_2 + \left[  v_{23}(t) \varphi_{2,m}(l_2) - v_{12}(t) \varphi_{2,m}(l_{1})\right],\\
\frac{\partial \widetilde{X}_3}{\partial t} &=-D_{b3} \lambda_{3,m}^2 \widetilde{\mathcal{X}}_3 + \left[ v_{13}(t) \varphi_{3,m}(l_{3}) -v_{23}(t) \varphi_{3,m}(l_2)  \right],
\end{align*}
together with the transformed initial conditions $\widetilde{X}_i(\lambda_{i,m},0) = \langle X_{i,0}(x), \varphi_{i,m} \rangle, i=1,2,3$. We now apply the Laplace transform in time and denote $\overline{\widetilde{X}}_i(\lambda_{i,m},s)=\mathcal{L}\left\{ \widetilde{X}_i(\lambda_{i,m},t) \right\},i=1,2,3$.  Since $\mathcal{L} \left\{ {_0^{R}D}_{t}^{1-\gamma_i} {\widetilde{X}}_i \right\} =s^{1-\gamma_i}  {\widetilde{X}}_i(\lambda_{i,m},s)$, the Laplace transformation of the three layer equations given above can be rearranged as
\begin{align}
\overline{\widetilde{X}}_1(\lambda_{1,m},s)=\frac{\widetilde{X}_1(\lambda_{1,m},0)}{\eta_{1,m}(s)} &+ \frac{ \bar{v}_{12}(s) \varphi_{1,m}(l_1)-\bar{v}_{13}(s) \varphi_{1,m}(l_0) }{\eta_{1,m}(s)}, \label{eq18} \\
\overline{\widetilde{X}}_2(\lambda_{2,m},s)=\frac{\widetilde{X}_2(\lambda_{2,m},0)}{\eta_{2,m}(s)} &+ \frac{ \bar{v}_{23}(s) \varphi_{2,m}(l_2) - \bar{v}_{12}(s) \varphi_{2,m}(l_1) }{\eta_{2,m}(s)}, \label{eq19} \\
\overline{\widetilde{X}}_3(\lambda_{3,m},s)=\frac{\widetilde{X}_3(\lambda_{3,m},0)}{\eta_{3,m}(s)} &+ \frac{  \bar{v}_{13}(s) \varphi_{3,m}(l_3)-\bar{v}_{23}(s) \varphi_{3,m}(l_2) }{\eta_{3,m}(s)}, \label{eq20}
\end{align}
where $\eta_{i,m}(s)=s+D_{bi} s^{1-\gamma_i} \lambda_{i,m}^2$. In order to determine the three unknown interfacial flux values $\bar{v}_{12}(s)$, $\bar{v}_{13}(s)$ and $\bar{v}_{23}(s)$, we need the boundary conditions at the interfaces: $\overline{X}_1(l_1,s)=\overline{X}_2(l_1,s)$, $\overline{X}_3(l_3,s)=\overline{X}_1(l_0,s)+\frac{q_0}{s}$ and $\overline{X}_2(l_2,s)=\overline{X}_3(l_2,s)$. Noting that $\overline{X}_i(x,s)=\sum_{m=0}^{\infty}\overline{\widetilde{X}}_i(\lambda_{i,m},s) \varphi_{i,m}(x)$, $i=1,2,3$ and substituting the expressions (\ref{eq18})-(\ref{eq20}) we obtain the following $3 \times 3$ linear system:
\begin{align}
\sum_{m=0}^{\infty}\overline{\widetilde{X}}_1(\lambda_{1,m},s) \varphi_{1,m}(l_1)&= \sum_{m=0}^{\infty}\overline{\widetilde{X}}_2(\lambda_{2,m},s) \varphi_{2,m}(l_1), \nonumber\\
\sum_{m=0}^{\infty}\overline{\widetilde{X}}_3(\lambda_{3,m},s) \varphi_{1,m}(l_3)&= \sum_{m=0}^{\infty}\overline{\widetilde{X}}_1(\lambda_{1,m},s) \varphi_{2,m}(l_0)+\frac{q_0}{s}, \nonumber\\\label{eq21}
\sum_{m=0}^{\infty}\overline{\widetilde{X}}_2(\lambda_{2,m},s) \varphi_{2,m}(l_2)&= \sum_{m=0}^{\infty}\overline{\widetilde{X}}_3(\lambda_{3,m},s) \varphi_{3,m}(l_2),
\end{align}
that can be solved for $\bar{v}_{12}(s)$, $\bar{v}_{13}(s)$ and $\bar{v}_{23}(s)$ at a given value of $s$ using Cramer's rule. Finally, the solutions within each layer can be determined by applying the inverse Laplace transform resolved numerically:
\begin{align*}
{\widetilde{X}}_i(\lambda_{i,m},t)&=\mathcal{L}^{-1}\left\{ {\overline{\widetilde{X}}}_i(\lambda_{i,m},s) \right\} = \frac{1}{2 \pi i} \int_{\Gamma} e^{st} \overline{{\widetilde{X}}}_i(\lambda_{i,m},s) \, ds\\
&= \frac{1}{2 \pi i} \int_{\Gamma} \frac{e^z}{t}  \overline{{\widetilde{X}}}_i(\lambda_{i,m},z/t) \, dz \approx -2\Re\left( \sum_{k=1}^{K/2}c_{2k-1} \frac{ \overline{{\widetilde{X}}}_i(\lambda_{i,m},z_{2k-1}/t)}{t}\right),
\end{align*}
where $z=st$, $c_{2k-1}$ and $z_{2k-1}$ are the residues and poles of the best $(K,K)$ rational approximation of $e^z$ on the negative real line as computed by the Carath\'{e}odory--Fej\'{e}r method. Full details of this numerical approach can be found in Trefethen et al. \cite{Tref06}. Then the solution in each layer can be computed using $X_i(x,t)=\sum_{m=0}^{\infty}\widetilde{X}(\lambda_{i,m},t) \varphi_{i,m}(x)$, $i=1,2,3$. We will use this semi-analytical solution to validate the numerical solutions we compute using the finite volume method summarised in Section \ref{sec2}.

\subsection{Layered medium mass balance equations}\label{sec3.3}

We now determine mass balance equations to confirm the accuracy of our simulation results. Define $L=\sum_{j=1}^3 d_j$ and the average value of the variable $X$ as follows:
\begin{equation*}
\langle X \rangle = \frac{1}{L} \int_{l_0}^{l_3} X(x,t) \, dx =  \frac{1}{L}  \sum_{i=1}^3 \int_{l_{i-1}}^{l_i} X_{i}(x,t) \, dx = \sum_{i=1}^3 \left( \frac{d_i}{L} \right) \langle X_i \rangle,
\end{equation*}
where the average value of $X$ in the $i^{\mathrm{th}}$ layer is given by $\langle X_i \rangle = \frac{1}{d_i} \int_{l_{i-1}}^{l_i} X_i(x,t) \, dx$. Using the continuity of the fluxes at the interfaces, we now integrate (\ref{eq16}) across each layer to obtain:
\begin{align} \label{eq22}
\frac{d \langle X \rangle}{dt} = 0,
\end{align}
subject to $\langle X \rangle(0) = \frac{1}{L} \int_{l_0}^{l_3} X_0(x) \, dx = \langle X_0 \rangle = \sum_{i=1}^3 \left( \frac{d_i}{L} \right) \langle X_{i0} \rangle$, where the $i^{\mathrm{th}}$ layer average initial condition is $ \langle X_{i0} \rangle = \frac{1}{d_i} \int_{l_{i-1}}^{l_i} X_{i0}(x,0) \, dx$. Then we can obtain that the average value of $X$ must remain constant at its initial average for the duration of the simulation, namely
\begin{align} \label{eq23}
 \langle X \rangle &=  \langle X_0 \rangle.
\end{align}
\section{Homogenisation theory for time-fractional differential equation}\label{sec4}

In this section, we will extend the homogenisation theory from classical differential equations to fractional differential equations in heterogeneous media. Let $\Omega$ be a bounded domain in $\mathbb{R}^2$ with (sufficiently smooth) boundary $\partial \Omega$. At first, we note that the following two equations are equivalent.
\begin{align} \label{eq24}
&\frac{\partial u}{\partial t}={_0^RD_t^{1-\gamma}} \left[\nabla(D\cdot\nabla u) \right],\\\label{eq25}
&{_0^CD_t^{\gamma}}u=\left[\nabla(D\cdot\nabla u) \right].
\end{align}
This equivalence can be verified by imposing the operator ${_0^RD_t^{\gamma-1}}$ on both sides of (\ref{eq24}) and using the relationship between the Riemann-Liouville fractional operator and the Caputo fractional operator \cite{Zhuang2} to obtain (\ref{eq25}). Due to this equivalence, for ease of mathematical exposition we focus on the homogenisation theory for the fractional differential equation (\ref{eq25}) involving the Caputo fractional derivative.  We consider the following generalised fractional transport model:
\begin{align}\label{eq26}
&{_0^CD_t^\gamma}u^\varepsilon(x,t)=\nabla\cdot(D^\varepsilon(x,t)\nabla u^\varepsilon(x,t)),~ \text{in}~\Omega,~t\in(0,T],\\ \nonumber
&u^\varepsilon(x,0)=f(x),~ \text{in}~\Omega,\quad u^\varepsilon(x,t)=0,~ \text{on}~\partial\Omega,~t\in(0,T],
\end{align}
where the Caputo derivative is defined as ${_0^CD_t^\gamma}u(x,t)=\frac{1}{\Gamma(1-\gamma)}\int_0^t (t-\eta)^{-\gamma}\frac{\partial u(x,\eta)}{\partial \eta}d\eta$, and $0<\gamma<1$ is the fixed time fractional derivative index. We assume the diffusivity tensor $D^\varepsilon(x,t)$ is symmetric with $D^\varepsilon_{ij}\in L^\infty(\Omega)$ and $D^\varepsilon_{ij}(x,t)\xi_i\xi_j\geq C\xi_i^2$ a.e. in $\Omega$, for some constant $C>0$, $\xi_i\in\mathbb{R}$. We define $D^\varepsilon(x,t):=D \left(\frac{x}{\varepsilon}, \frac{t}{\varepsilon^p}\right)$, $p>0$. The choice of $p$ will be discussed below and we assume there is no periodicity in the time variables. We follow the discussion from \cite{Bensoussan} by assuming further conditions on the structure of the functions $D^\varepsilon_{ij}$, namely:
\begin{itemize}
\item[1.] $D^\varepsilon_{ij}(x,t):=D_{ij} \left(\frac{x}{\varepsilon}, \frac{t}{\varepsilon^p}\right)$, $p>0$.
\item[2.] $D_{ij}(y,\tau)$ is $Y$-periodic as a function of $y$, where the fast variables in space and time are $y=\frac{x}{\varepsilon}$ and $\tau=\frac{t}{\varepsilon^p}$, respectively.
\item[3.] $D_{ij}(y,\tau)\in L^\infty(\mathbb{R}_y^2\times\mathbb{R}_\tau)$.
\end{itemize}
Assuming the fast variable $y$ and slow variable $x$ are independent as $\varepsilon\to 0$, we seek an asymptotic expansion of the solution $u^\varepsilon(x,t)$ in the form $u^\varepsilon(x,t)=\sum_{k=0}^{\infty} \varepsilon^ku_k(x,y,t,\tau)$, with $u_k$ being $Y$-periodic in $y$ and we will add further constraints on $u_k$ regarding the time variable $\tau$ below. Before proceeding with the two-scale asymptotic expansion of the solution $u^\varepsilon$, we first consider a general function $f(t,\tau)$, with $\tau=\frac{t}{\varepsilon^p}$, $p>0$ being the local scale in time. Actioning the Caputo derivative gives
\begin{align}\label{eq27}
{_0^C\mathbb{D}_t^\gamma}f(t,\tau)=\frac{1}{\Gamma(1-\gamma)}\int_0^t (t-\eta)^{-\gamma}f'(\eta,\tau)d\eta.
\end{align}
Next, the chain rule gives $\frac{d}{dt}=\frac{\partial}{\partial t}+\frac{1}{\varepsilon^p}\frac{\partial}{\partial\tau}$, i.e., $f'(\eta,\tau)=\frac{\partial f}{\partial t}+\frac{1}{\varepsilon^p}\frac{\partial f}{\partial\tau}$. Substituting this expression into (\ref{eq27}), we obtain
\begin{align*}
{_0^C\mathbb{D}_t^\gamma}f(t,\tau)&=\frac{1}{\Gamma(1-\gamma)}\int_0^t (t-\eta)^{-\gamma} \left(\frac{\partial f(\eta,\zeta)}{\partial t}+\frac{1}{\varepsilon^p}\frac{\partial f(\eta,\zeta)}{\partial\zeta} \right)d\eta\\
&={_0^CD_t^\gamma}f(t,\tau)+\frac{1}{\Gamma(1-\gamma)\varepsilon^p}\int_0^t (t-\eta)^{-\gamma}\frac{\partial f(\eta,\zeta)}{\partial\zeta}d\eta.
\end{align*}
The second term on the RHS can be written as
\begin{align*}
\frac{1}{\Gamma(1-\gamma)\varepsilon^p}\int_0^{\frac{t}{\varepsilon^p}} (t-\varepsilon^p\zeta)^{-\gamma}\frac{\partial f(\eta,\zeta)}{\partial\zeta}\varepsilon^pd\zeta=\frac{1}{\Gamma(1-\gamma)}\int_0^\tau\frac{\partial f(\eta,\zeta)}{\partial\zeta}\frac{d\zeta}{(\tau-\zeta)^{\gamma}\varepsilon^{p\gamma}},
\end{align*}
where we have set $\zeta=\frac{\eta}{\varepsilon^p}$ and $t=\varepsilon^p\tau$. Our investigation of different choices for the parameter $p$ identified that $p=\frac{2}{\gamma} $ gives
${_0^C\mathbb{D}_t^\gamma}f(t,\tau)={_0^CD_t^\gamma}f(t,\tau)+\frac{1}{\varepsilon^2}{_0^CD_\tau^\gamma}f(t,\tau)$. We will use this relation below to derive the time-fractional unit cell model. \medskip

\noindent Next, denote $A^\varepsilon=-\frac{\partial}{\partial x_i}\left[D_{ij}(y,\tau) \frac{\partial}{\partial x_j}  \right]$ and then expand $A^\varepsilon$ as:
\begin{align*}
A^\varepsilon=\varepsilon^{-2}A_1+\varepsilon^{-1}A_2+\varepsilon^{0}A_3,
\end{align*}
where $A_1=-\frac{\partial}{\partial y_i}\left[D_{ij}(y,\tau) \frac{\partial}{\partial y_j}  \right]$, $A_2=-\frac{\partial}{\partial y_i}\left[D_{ij}(y,\tau) \frac{\partial}{\partial x_j}  \right]-\frac{\partial}{\partial x_i}\left[D_{ij}(y,\tau) \frac{\partial}{\partial y_j}  \right]$, $A_3=-\frac{\partial}{\partial x_i}\left[D_{ij}(y,\tau) \frac{\partial}{\partial x_j}  \right]$. Substituting the expressions for $u^\varepsilon$ and $A^\varepsilon$ into system (\ref{eq24}) gives:
\begin{align*}
&\text{LHS:}~\left( {_0^CD_t^\gamma}+\frac{1}{\varepsilon^2}{_0^CD_\tau^\gamma}\right)(u_0+\varepsilon u_1+\varepsilon^2 u_2+\ldots)=\varepsilon^{-2}{_0^CD_\tau^\gamma}u_0+\varepsilon^{-1}{_0^CD_\tau^\gamma}u_1+\varepsilon^0\left( {_0^CD_t^\gamma}u_0+{_0^CD_t^\gamma}u_2 \right) +\ldots\\
&\text{RHS:}~-A^\varepsilon(u_0+\varepsilon u_1+\varepsilon^2 u_2+\ldots)=-\varepsilon^{-2}A_1u_0-\varepsilon^{-1}(A_1u_1+A_2u_0)-\varepsilon^0(A_1u_2+A_2u_1+A_3u_0)+\ldots
\end{align*}
Then by equating powers of $\varepsilon$, we obtain:
\begin{align}\label{eq28}
&O\left(\varepsilon^{-2} \right):~{_0^CD_\tau^\gamma}u_0=-A_1u_0, \\ \label{eq29}
&O\left(\varepsilon^{-1} \right):~{_0^CD_\tau^\gamma}u_1=-(A_1u_1+A_2u_0), \\ \label{eq30}
&O\left(\varepsilon^{0} \right):~{_0^CD_t^\gamma}u_0+{_0^CD_\tau^\gamma}u_2=-(A_1u_2+A_2u_1+A_3u_0).
\end{align}
Due to the complexity introduced by the time-fractional operator and the time dependent periodic boundary conditions at the external and interfacial boundaries, we choose to solve these problems to steady-state in $\tau$  at which point $u_0=u(x,t)$ is independent of $y$ and $\tau$ and (\ref{eq29}) can be expressed as
\begin{align*}
{_0^CD_\tau^\gamma}u_1+A_1u_1=\frac{\partial D_{ij}}{\partial y_i}\frac{\partial u}{\partial x_j}.
\end{align*}
Next, we introduce $\theta_j$ as the steady state $Y$-periodic solution of the following time-fractional problem
\begin{align*}
{_0^CD_\tau^\gamma}\theta_j=-A_1(\theta_j+y_j),\quad j=1,2.
\end{align*}
The solution $u_1$ can be defined up to an additive constant as
\begin{align*}
u_1(x,y,t,\tau)=\theta_j\frac{\partial u}{\partial x_j}+\overline{u}_1(x,t).
\end{align*}
Hence, the unit cell problem becomes to solve:
\begin{align*}
{_0^CD_\tau^\gamma}\theta_j=\frac{\partial}{\partial y_i}\left[D_{ij}(y,\tau) \frac{\partial (\theta_j+y_j)}{\partial y_j}  \right],\quad j=1,2,
\end{align*}
subject to periodic boundary conditions in $y$ and initially $\theta_j(y,0)$ is specified. Denote by $I_Y(\cdot)=\frac{1}{|Y|}\int_Y \cdot dy $ the average of a quantity in $y$. At steady-state in $\tau$ (\ref{eq30}) becomes
\begin{align*}
{_0^CD_t^\gamma}u=-I_y(A_2u_1+A_3u),
\end{align*}
which gives the homogenised equation as
\begin{align*}
{_0^CD_t^\gamma}u=\frac{\partial}{\partial x_i}\left(\overline{{D}}_{ij}(y,\tau) \frac{\partial u}{\partial x_j}  \right),~{\rm{where}}~ {\bm{\overline{D}}}=\frac{1}{|Y|}\int_Y {\bm{D}}({\bm{I}}+{\bm{J}}_\theta^T) dy,
\end{align*}
\begin{align*}
\bm{D}=\left[
         \begin{array}{cc}
           D_{11} & D_{12} \\
           D_{21} & D_{22} \\
         \end{array}
       \right],\quad
\bm{J}_\theta = \left[
         \begin{array}{cc}
           \frac{\partial\theta_1}{\partial y_1} & \frac{\partial\theta_1}{\partial y_2} \\
           \frac{\partial\theta_2}{\partial y_1} & \frac{\partial\theta_2}{\partial y_2} \\
         \end{array}
       \right].
\end{align*}
We can express the unit cell problem in terms of the new variable $\varphi_j=\theta_j+y_j$, $j=1,2$, as
\begin{align}\label{eq31}
{_0^CD_\tau^\gamma}\varphi_j=\frac{\partial}{\partial y_i}\left[D_{ij}(y,\tau) \frac{\partial \varphi_j}{\partial y_j}  \right].
\end{align}
This variable change affects the boundary conditions:
\begin{align*}
\varphi_1(L_1,y_2,\tau)&=\varphi_1(0,y_2,\tau)+L_1,\quad \varphi_1(y_1,L_2,\tau)=\varphi_1(y_1,0,\tau),\\
\varphi_2(L_1,y_2,\tau)&=\varphi_2(0,y_2,\tau),\quad \varphi_2(y_1,L_2,\tau)=\varphi_2(y_1,0,\tau)+L_2.
\end{align*}
Also we need the boundary conditions (\ref{eq12}) at the interface. The initial conditions become $\varphi_j(y_1,y_2,0)=\theta_j(y_1,y_2,0)+y_j$, $j=1,2$. In this case we obtain the effective diffusivity tensor as
\begin{align*}
{\bm{\overline{D}}}=\frac{1}{|Y|}\int_Y {\bm{D}}\bm{J}_\varphi^Tdy.
\end{align*}

\section{Numerical examples}\label{sec5}
\begin{figure}[H]
\subfloat[Morphology 1: Rectangular inclusion]{
\label{Fig7a}
\begin{minipage}[t]{1\textwidth}
\centering{
\scalebox{0.3}[0.3]{\includegraphics{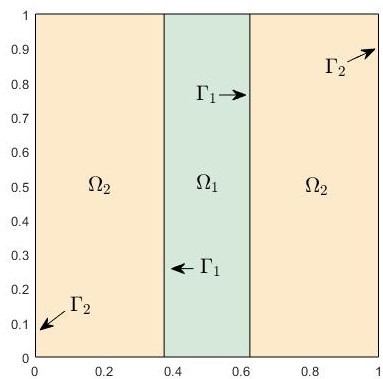}}~
\scalebox{0.3}[0.3]{\includegraphics{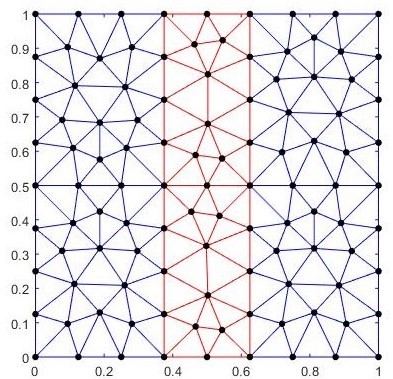}}~
\scalebox{0.3}[0.3]{\includegraphics{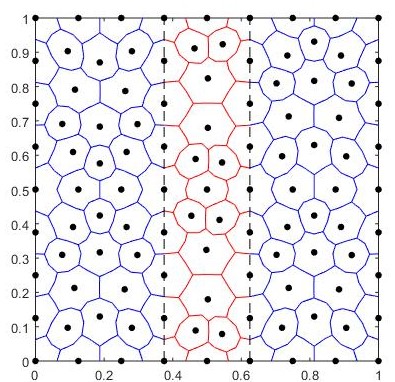}}}
\end{minipage}}\\
\subfloat[Morphology 2: Circular inclusion]{
\label{Fig7b}
\begin{minipage}[t]{1\textwidth}
\centering{
\scalebox{0.315}[0.315]{\includegraphics{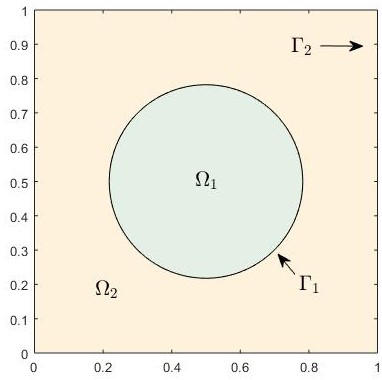}}~
\scalebox{0.3}[0.3]{\includegraphics{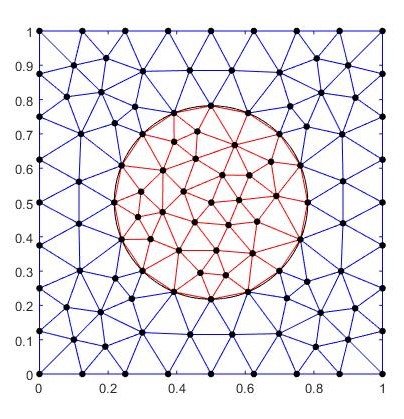}}~
\scalebox{0.3}[0.3]{\includegraphics{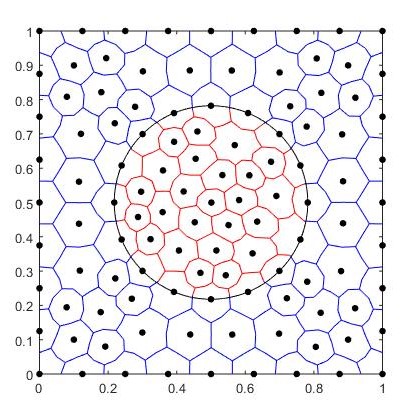}}}
\end{minipage}}\\
\subfloat[Morphology 3: L-shaped inclusion]{
\label{Fig7c}
\begin{minipage}[t]{1\textwidth}
\centering{
\scalebox{0.3}[0.3]{\includegraphics{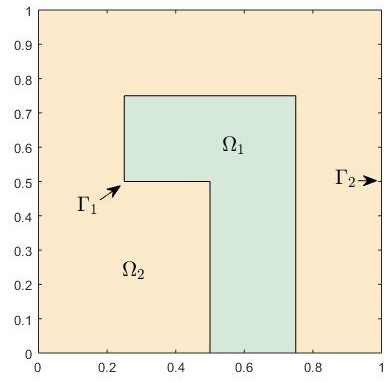}}~
\scalebox{0.31}[0.31]{\includegraphics{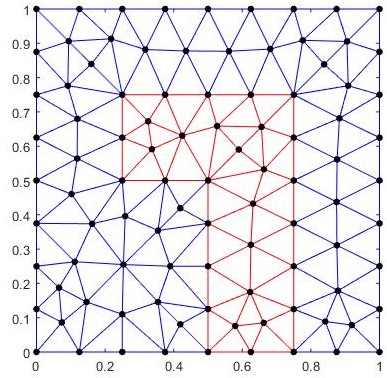}}~
\scalebox{0.31}[0.31]{\includegraphics{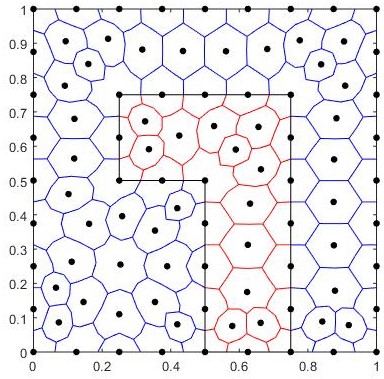}}}
\end{minipage}}
\caption{Schematics of three binary media taken from \cite{Szymk}, each with volumetric fraction for $\Omega_1$ taken as 0.25. The triangulation and control volume partitions are also shown. (a) Morphology 1 $\Omega_1=[\frac{3}{8},\frac{5}{8}]\times [0,1]$, $\Omega_2=[0,1]\times[0,1]\backslash \Omega_1$. (b) Morphology 2 $\Omega_1=\{(x,y)|(x-0.5)^2+(y-0.5)^2<\frac{1}{4\pi}\}$, $\Omega_2=[0,1]\times[0,1]\backslash \Omega_1$. (c) Morphology 3 $\Omega_1=[\frac{1}{2},\frac{3}{4}]\times [0,\frac{1}{2}]\cup [\frac{1}{4},\frac{3}{4}]\times [\frac{1}{2},\frac{3}{4}]$ , $\Omega_2=[0,1]\times[0,1]\backslash \Omega_1$.}
\label{Fig7}
\end{figure}

In this section, we present the results obtained from the computational homogenisation simulations for binary media having the three different morphologies shown in Figure \ref{Fig7}: Morphology 1 is a two-layered problem with a rectangular inclusion; Morphology 2 is a binary medium with a circular inclusion; and Morphology 3 is a binary medium with an L-shaped inclusion. All of the computations were carried out using MATLAB R2018a on a DELL desktop with the configuration: Intel(R) Core(TM) i7-6700 CPU\@3.40GHz and RAM 16.0 GB.

To perform the homogenisation, we solve (\ref{eq31}) with periodic boundary conditions and nonzero initial conditions. The effective diffusivity (\ref{eq32}) is computed once steady-state is achieved. For the homogenisation problems considered we take $D_{12}=D_{21}=0$ and denote
\begin{align}\label{eq32}
{\bm{\overline{D}}}= \left[
                                \begin{array}{cc}
                                  D_{b_{x}}  &  D_{b_{xy}} \\
                                  D_{b_{yx}} &  D_{b_{y}} \\
                                \end{array}
                              \right].
\end{align}

\noindent In order to ensure the correct interfacial boundary conditions are invoked, we revert to the following time-fractional model with the Riemann-Liouville fractional derivative:
\begin{equation*}
\left\{\begin{array}{l}
\frac{\partial u(x,y,t)}{\partial t}={_0^RD^{1-\gamma_1}_{t}}\left[\nabla\cdot(D_{b_1}\nabla u(x,y,t))\right],\quad(x,y,t)\in \Omega_1\times (0,T],\\
\frac{\partial v(x,y,t)}{\partial t}={_0^RD^{1-\gamma_2}_{t}}\left[\nabla\cdot(D_{b_2}\nabla v(x,y,t))\right],\quad(x,y,t)\in \Omega_2\times (0,T],
\end{array}\right.
\end{equation*}
with Neumann boundary conditions
\begin{align*}
{_0^RD^{1-\gamma_1}_{t}}\left[-D_{b_1}\nabla u(x,0,t)\cdot {\bf{n}}\right]&={_0^RD^{1-\gamma_1}_{t}}\left[D_{b_1}\nabla u(x,1,t)\cdot {\bf{n}}\right],\\
{_0^RD^{1-\gamma_2}_{t}}\left[-D_{b_2}\nabla v(0,y,t)\cdot {\bf{n}}\right]&={_0^RD^{1-\gamma_2}_{t}}\left[D_{b_2}\nabla v(1,y,t)\cdot {\bf{n}}\right],\\
{_0^RD^{1-\gamma_2}_{t}}\left[-D_{b_2}\nabla v(x,0,t)\cdot {\bf{n}}\right]&={_0^RD^{1-\gamma_2}_{t}}\left[D_{b_2}\nabla v(x,1,t)\cdot {\bf{n}}\right],
\end{align*}
and quasi-periodic boundary conditions $v(1,y,t)=v(0,y,t)+1,~ v(x,1,t)=v(x,0,t)+1$, and with additional boundary conditions at the interface $\Gamma_1$ are $ u(x,y,t)= v(x,y,t)$, ${_0^RD^{1-\gamma_1}_{t}}(D_{b_1}\nabla u(x,y,t)\cdot {\bf{n}})={_0^RD^{1-\gamma_2}_{t}}(D_{b_2}\nabla v(x,y,t)\cdot {\bf{n}})$.\medskip

We use the time evolutionary homogenisation theory from Section \ref{sec4} to estimate the effective parameters for these three different morphologies as comparisons with the values of the diffusivity tensors can be made with the test problems given in \cite{Szymk}. To calculate the equivalent diffusivity, we advance the numerical scheme in time until the steady-state is reached. We then calculate the equivalent diffusivity $\overline{D_{b_x}}$ in the $x$-direction and $\overline{D_{b_y}}$ in the $y$-direction using the following formulae obtained by approximating (\ref{eq32}) in Section \ref{sec4} using a midpoint quadrature rule applied to the integral over the unit cell
\begin{align*}
\overline{D_{b_x}}=\frac{\sum\limits_{i}\left(\chi D_{b_1}\frac{\partial u}{\partial x}+(1-\chi)D_{b_2}\frac{\partial v}{\partial x}\right)_iS_{\triangle_i}}{\sum\limits_{i}S_{\triangle_i}},\quad
\overline{D_{b_y}}=\frac{\sum\limits_{i}\left(\chi D_{b_1}\frac{\partial u}{\partial y}+(1-\chi)D_{b_2}\frac{\partial v}{\partial y}\right)_iS_{\triangle_i}}{\sum\limits_{i}S_{\triangle_i}},
\end{align*}
where $S_{\triangle_i}$ is the area of the $i^{\mathrm{th}}$ triangular element and the indicator function is defined as
\begin{equation*}
\chi=\left\{\begin{array}{ll}
1,& (x,y)\in \Omega_1,\\
0,& (x,y)\in \Omega_2.
\end{array}\right.
\end{equation*}
Throughout this section we will frequently make reference to the harmonic average ($K_1$) and arithmetic average ($K_2$), which can be calculated by $K_1=\left(\frac{\epsilon_1}{D_{b_1}}+\frac{\epsilon_2}{D_{b_2}}\right)^{-1}$, $K_2=\epsilon_1D_{b_1}+\epsilon_2D_{b_2}$, where $\epsilon_1$ is the volumetric fraction of $\Omega_1$ and  $\epsilon_2$ is the volumetric fraction of $\Omega_2$. \medskip

\noindent {\bf{Morphology 1}}: \emph{Rectangular inclusion (layered composite material)}.\medskip

\begin{figure}[H]
\centering
\subfloat[Different $\gamma_1$]{
\label{Fig4a}
\begin{minipage}[t]{0.33\textwidth}
\centering
\scalebox{0.4}[0.4]{\includegraphics{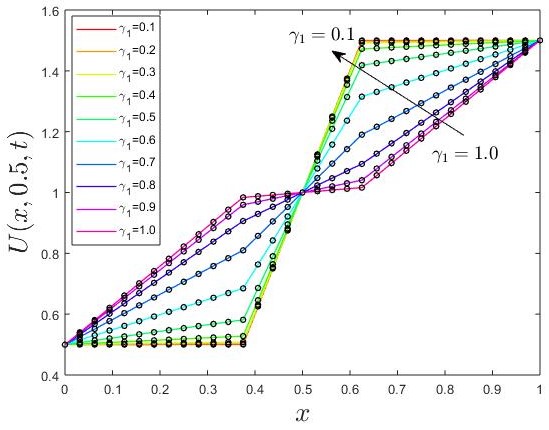}}
\end{minipage}
}
\subfloat[Different $t$]{
\label{Fig4b}
\begin{minipage}[t]{0.33\textwidth}
\centering
\scalebox{0.43}[0.43]{\includegraphics{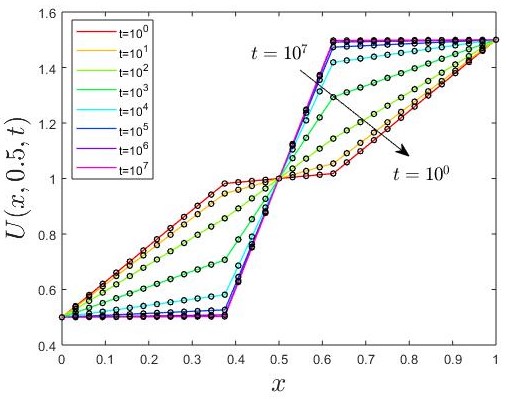}}
\end{minipage}
}
\subfloat[Different $\gamma_2$]{
\label{Fig4c}
\begin{minipage}[t]{0.33\textwidth}
\centering
\scalebox{0.43}[0.43]{\includegraphics{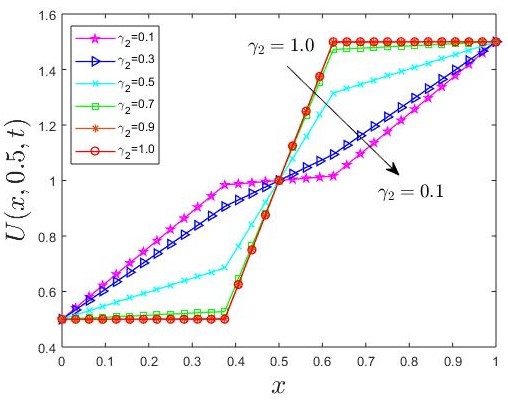}}
\end{minipage}
}
\caption{A comparison between the numerical solution (symbol) and semi-analytical solution (line) for different $\gamma_1$, $\gamma_2$ and $t$ with $D_{b_1}=10$, $D_{b_2}=1$ and $q_0=u_0=1$, in which other parameters are: (a) $\gamma_2=1$ at $t=10^4\,s$; (b) $\gamma_1=0.5$ and $\gamma_2=1$; (c) $\gamma_1=0.1$ at $t=10^4\,s$.}
\label{Fig4}
\end{figure}

\noindent For Morphology 1, we are able to use the derived semi-analytic solution from Section \ref{sec3.2} to validate the CVM method to simulate time-fractional transport in a layered medium with quasi-periodic boundary conditions. This problem is chosen to validate the CVM method when applied to the homogenisation theory presented in Section \ref{sec4}. We exhibit the numerical solution at the central line $y=0.5$ to compare directly with the semi-analytical solution derived in Section \ref{sec3.2}. The initial condition is taken as $u(x,y,0)=v(x,y,0)=u_0$. The maximum error between the numerical and semi-analytical solutions is presented in Table \ref{tab2} for different $\gamma_1$ and $\gamma_2$ with $h=1.848\times 10^{-1}$, $\tau=1\times10^{-3}$, $D_{b_1}=10$, $D_{b_2}=1$, $u_0=q_0=1$ at $t=1$. We conclude from the results in Table \ref{tab2} and Figure \ref{Fig4} that the numerical and semi-analytical solutions agree very well, with the maximum error in all cases $\approx 3.2 \times 10^{-4}$.
\begin{table}[H]
\begin{center}
\caption{The maximum error between the numerical solution and semi-analytical solution at the central line $y=0.5$ for different $\gamma_1$ and $\gamma_2$ with $h=1.848\times 10^{-1}$, $\tau=1\times10^{-3}$, $D_{b_1}=10$, $D_{b_2}=1$, $u_0=q_0=1$ at $t=1$.}
\label{tab2}
\begin{tabular}{cccc}
\toprule
      &  $\gamma_1$=0.2 & $\gamma_1$=0.5 & $\gamma_1$=0.8\\
\midrule
$\gamma_2$=0.2  & 3.2198E-04  & 3.2206E-04 & 3.2088E-04   \\
$\gamma_2$=0.5  & 3.2002E-04  & 3.1897E-04 & 3.1818E-04  \\
$\gamma_2$=0.8  & 3.2182E-04  & 3.1903E-04 & 3.1749E-04  \\
$\gamma_2$=1.0  & 3.2632E-04  & 3.2194E-04 & 3.1917E-04  \\
\bottomrule
\end{tabular}
\end{center}
\end{table}
\vspace{-5mm}
Next, we investigate the impact of varying the fractional indices $\gamma_1$ and $\gamma_2$ on the solution profile at steady-state. At first, we fix $\gamma_2=1$ to reduce the diffusion on $\Omega_2$ to normal diffusion and then change $\gamma_1$ to observe the solution behaviour at $t=10^4\,s$ (see Figure \ref{Fig4a}). At this time, the diffusion for $\gamma_1=1$ has reached its steady-state, which corresponds with the solution of the classical two-layered diffusion problem. The diffusion for $\gamma_1=0.1$ also reaches its steady-state, however the solution profile is totally different to the classical case, being much steeper at the centre point $x=0.5$.
\vspace{-3mm}
We conclude that the fractional index $\gamma_1$ has a significant impact on the diffusion process for the two-layered media with memory. Another interesting finding is that a larger value of the fractional index $\gamma_1$ will delay the diffusion process from reaching steady-state. Furthermore, the final steady-state profiles for all the fractional cases are identical when $\gamma_2=1$. To investigate this phenomenon further, we now fix $\gamma_2=1$ and choose $\gamma_1=0.5$ to observe the solution behaviour from $t=1\,s$ to $t=10^7\,s$ (see Figure \ref{Fig4b}). We can see that the diffusion for $\gamma_1=0.5$ needs $t=10^7\,s$ to reach its steady-state, which is much longer than required for the case $\gamma_1=0.1$ ($t=10^4\,s$).  We again find that the steady-state solution profile for $\gamma_1=0.5$ is the same as exhibited in Figure \ref{Fig4a} for $\gamma_1=0.1$. Moreover, when we fix $\gamma_1=0.1$ to see the impact of $\gamma_2$ on the diffusion, the reverse phenomenon to that observed above is apparent (see Figure \ref{Fig4c}). It appears that the larger the choice of $\gamma_2$, the quicker it reaches the steady-state solution.
\begin{figure}[H]
\centering
\subfloat[Different $D_{b_1}$]{
\label{Fig5a}
\begin{minipage}[t]{0.33\textwidth}
\centering
\scalebox{0.43}[0.43]{\includegraphics{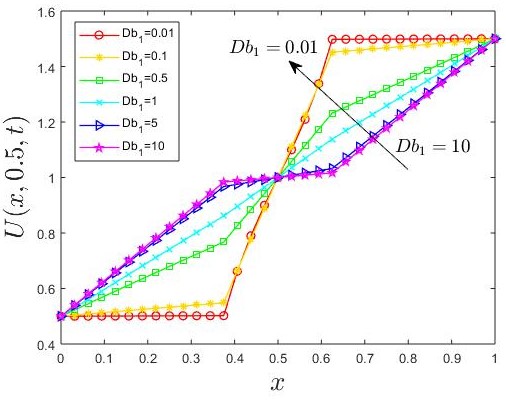}}
\end{minipage}
}
\subfloat[Different $q_0$]{
\label{Fig5b}
\begin{minipage}[t]{0.33\textwidth}
\centering
\scalebox{0.43}[0.43]{\includegraphics{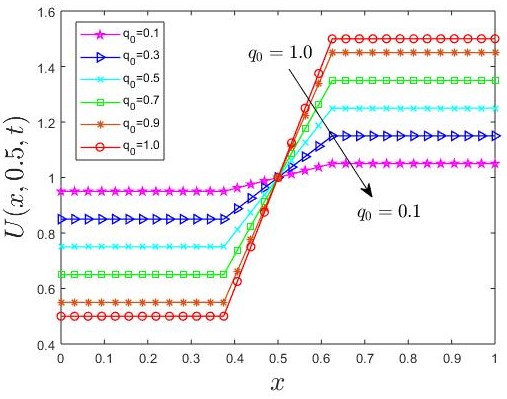}}
\end{minipage}
}
\subfloat[Different $u_0$]{
\label{Fig5c}
\begin{minipage}[t]{0.33\textwidth}
\centering
\scalebox{0.43}[0.43]{\includegraphics{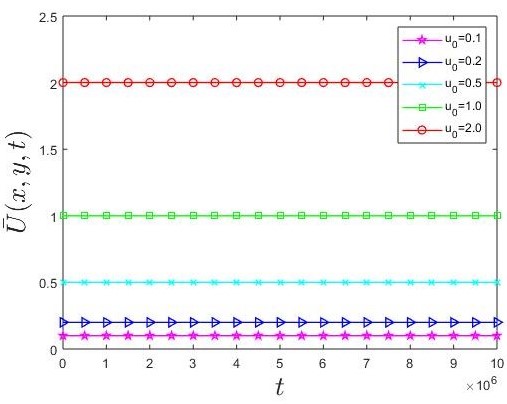}}
\end{minipage}
}
\caption{A comparison between the numerical solution (symbol) and semi-analytical solution (line) for different $D_{b_1}$, $q_0$ and $u_0$ with $\gamma_1=0.1$ and $\gamma_2=1$, in which other parameters are: (a) $D_{b_2}=1$, $q_0=u_0=1$ at $t=1\,s$; (b) $D_{b_1}=10$, $D_{b_2}=1$, $u_0=1$ at $t=10^4\,s$; (c) $D_{b_1}=10$, $D_{b_2}=1$, $q_0=1$ and $t$ from $0\,s$ to $10^7\,s$.}
\label{Fig5}
\end{figure}
\vspace{-3mm}
In Figure \ref{Fig5}, the effect of the other model parameters $D_{b1}, D_{b2}, q_0$ and $u_0$ on the diffusion process is exhibited. Figure \ref{Fig5a} shows the impact of the diffusivity coefficient for the memory part ($D_{b_1}$) with fixed $D_{b_2}=1$, $\gamma_1=0.1$, $\gamma_2=1$, from which we see that the smaller the value of $D_{b_1}$, the easier it becomes for the process to attain steady state. Figures \ref{Fig5b} and \ref{Fig5c} show the solution profiles for different values of $q_0$ and $u_0$ for both the numerical and semi-analytical solutions, which validates the accuracy of computational model and the mass balance equation (\ref{eq23}). In particular, it is evident that as the value of $q_0$ is decreased, the steady-state solution profile becomes flatter. Furthermore, for all choices of $u_0$, the average solution value remains constant at $u_0$ during the simulation.

We now give a comparison between the Riemann-Liouville (RL) and Caputo variants of the model discussed in Section \ref{sec4}:
\begin{equation}\label{eq35}
\left\{\begin{array}{l}
{_0^CD^{\gamma_1}_{t}} u(x,y,t)=D_{b_1}\Delta u(x,y,t),\quad(x,y,t)\in \Omega_1\times (0,T],\\
{_0^CD^{\gamma_2}_{t}} v(x,y,t)=D_{b_2}\Delta v(x,y,t),\quad(x,y,t)\in \Omega_2\times (0,T],
\end{array}\right.
\end{equation}
\noindent with the initial and periodic boundary conditions for the Caputo model now expressed in terms of the classical flux
\begin{align*}
u(x,y,0)=v(x,y,0)=u_0,~ D_{b_2}\frac{\partial v(0,y,t)}{\partial x}= D_{b_2}\frac{\partial v(1,y,t)}{\partial x},\quad v(1,y,t)= v(0,y,t)+q_0,
\end{align*}
and the boundary conditions at the interface $\Gamma_1$ between the two media
\begin{align}\label{eq36}
 u(x,y,t)= v(x,y,t), \quad D_{b_1}\frac{\partial u(x,y,t)}{\partial x}=D_{b_2}\frac{\partial v(x,y,t)}{\partial x},\quad (x,y,t)\in \Gamma_1\times (0,T].
\end{align}
\begin{figure}[H]
\begin{center}
\scalebox{0.35}[0.35]{\includegraphics{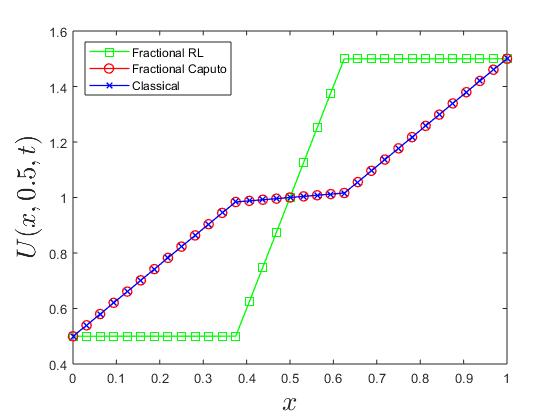}}
\caption{A comparison between the numerical solution (symbol) and semi-analytical solution (line) between the classical and fractional two-layered problems for the different time-fractional transport models, in which the the parameters are $\gamma_1=1$, $\gamma_2=1$ for the classical case and $\gamma_1=0.1$, $\gamma_2=1$ for the fractional case with $D_{b_1}=10$, $D_{b_2}=1$, $q_0=u_0=1$ at $t=10^4$.}
\label{Fig6}
\end{center}
\end{figure}
\vspace{-6mm}
As stated in Section \ref{sec4}, the fractional partial differential equations (\ref{eq24}) and (\ref{eq25}) are equivalent, however a most important observation is that the interfacial boundary conditions (\ref{eq12}) and (\ref{eq36}) are not. Following a similar approach as outlined in Section \ref{sec3.2}, we have derived the semi-analytical solution for the Caputo model (\ref{eq35}) and also exhibited its behaviour along the central line $y = 0.5$. Figure \ref{Fig6} compares the steady-state solution profiles between the classical and fractional two-layered problems for the two different time-fractional transport models. A key finding is that for the Caputo model (\ref{eq35}), there is no observable difference in the steady-state profiles computed for the classical and fractional cases whatever the choice of $\gamma_1$. This means that the impact of the internal layer memory effect is diminished, which may not be a true accord of the physics associated with the anomalous transport phenomena. The main reason for this outcome is due to the incorrect treatment of the interfacial boundary conditions (\ref{eq36}), which also explains why we need to transform the Caputo fractional derivative to the Riemann-Liouville fractional derivative to perform the computational homogenisation. Otherwise, the interfacial boundary conditions (\ref{eq12}) would need to be imposed on the Caputo model (\ref{eq35}), which is more complicated to calculate.

We complete this section with a summary of our findings for the homogenised diffusivities computed using the time-fractional homogenisation theory for the layered binary medium. To ease the mathematical exposition, we denote $\alpha_1=1-\gamma_1$ and $\alpha_2=1-\gamma_2$. At first we consider the case $\alpha_1=\alpha_2=0$, which reduces the problem to the classical two-layered problem. Figures \ref{Fig8a} and \ref{Fig8b} show the equivalent diffusivity, the harmonic average and arithmetic average curves of the two-layered medium with different diffusivity ratios $\frac{D_{b_1}}{D_{b_2}}$. We can see that $\overline{D_{b_x}}$ and $\overline{D_{b_y}}$ coincide with the harmonic average and arithmetic average, respectively, which is a well known result. When $\frac{D_{b_1}}{D_{b_2}}=10$, $\overline{D_{b_x}}=1.2903$, $\overline{D_{b_y}}=3.250$; and for $\frac{D_{b_1}}{D_{b_2}}=0.1$, $\overline{D_{b_x}}=0.3077$, $\overline{D_{b_y}}=0.775$. Table \ref{tab3} shows the relative error of our numerical results with the theoretical equivalent diffusivity in \cite{Szymk}, which illustrates the effectiveness of the proposed numerical scheme.

\begin{figure}[H]
\centering
\subfloat[Morphology 1: $\frac{D_{b_1}}{D_{b_2}}=10$]{
\label{Fig9a}
\begin{minipage}[t]{0.5\textwidth}
\centering
\scalebox{0.3}[0.3]{\includegraphics{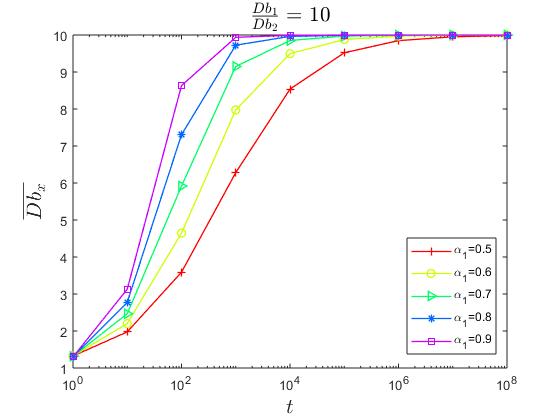}}
\end{minipage}}
\subfloat[Morphology 1: $\frac{D_{b_1}}{D_{b_2}}=0.1$]{
\label{Fig9b}
\begin{minipage}[t]{0.5\textwidth}
\centering
\scalebox{0.3}[0.3]{\includegraphics{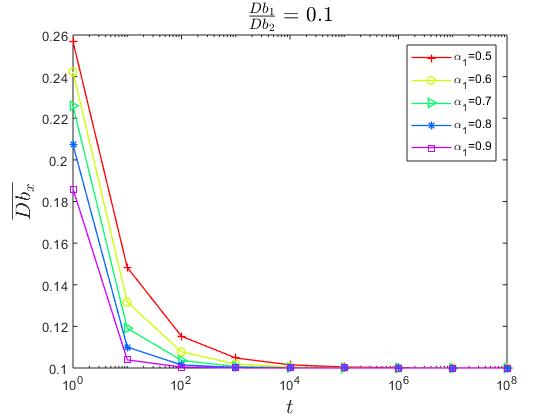}}
\end{minipage}}\\
\subfloat[Morphology 2: $\frac{D_{b_1}}{D_{b_2}}=10$]{
\label{Fig9c}
\begin{minipage}[t]{0.5\textwidth}
\centering
\scalebox{0.3}[0.3]{\includegraphics{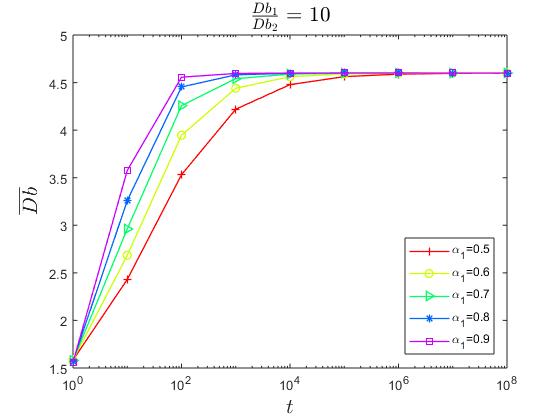}}
\end{minipage}}
\subfloat[Morphology 2: $\frac{D_{b_1}}{D_{b_2}}=0.1$]{
\label{Fig9d}
\begin{minipage}[t]{0.5\textwidth}
\centering
\scalebox{0.3}[0.3]{\includegraphics{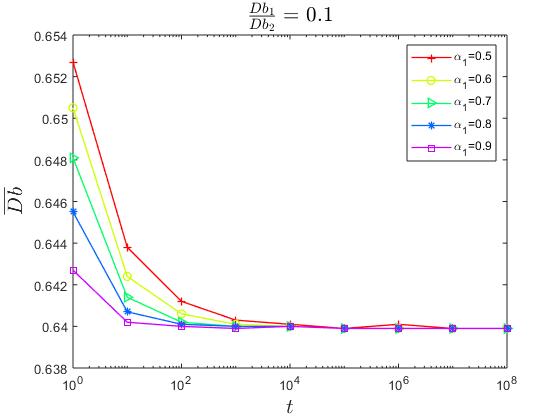}}
\end{minipage}}
\caption{The evolution of the equivalent diffusivity $\overline{D_{b_x}}$ using different $\alpha_1$ for Morphologies 1 and 2 with increasing time under different diffusivity ratios $\frac{D_{b_1}}{D_{b_2}}$, where $D_{b_2}=1$ and $\alpha_2=0$.}
\label{Fig9}
\end{figure}

\begin{figure}[H]
\centering
\subfloat[Morphology 1: $\frac{D_{b_1}}{D_{b_2}}=10$]{
\label{Fig8a}
\begin{minipage}[t]{0.35\textwidth}
\centering
\scalebox{0.35}[0.35]{\includegraphics{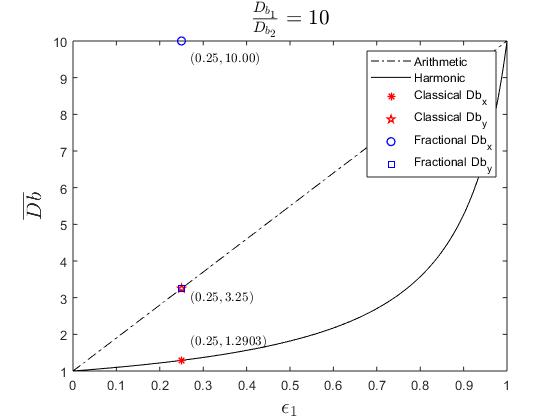}}
\end{minipage}}
\subfloat[Morphology 1: $\frac{D_{b_1}}{D_{b_2}}=0.1$]{
\label{Fig8b}
\begin{minipage}[t]{0.35\textwidth}
\centering
\scalebox{0.35}[0.35]{\includegraphics{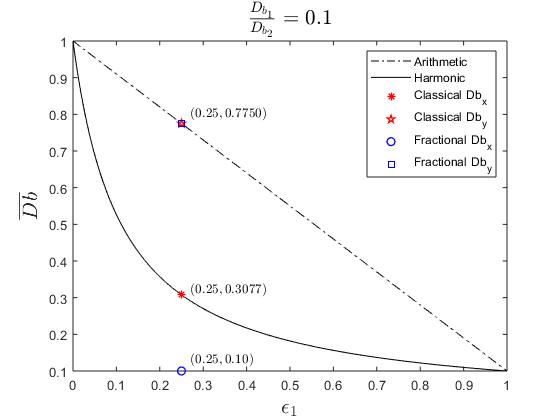}}
\end{minipage}}\\
\subfloat[Morphology 2: $\frac{D_{b_1}}{D_{b_2}}=10$]{
\label{Fig8c}
\begin{minipage}[t]{0.35\textwidth}
\centering
\scalebox{0.35}[0.35]{\includegraphics{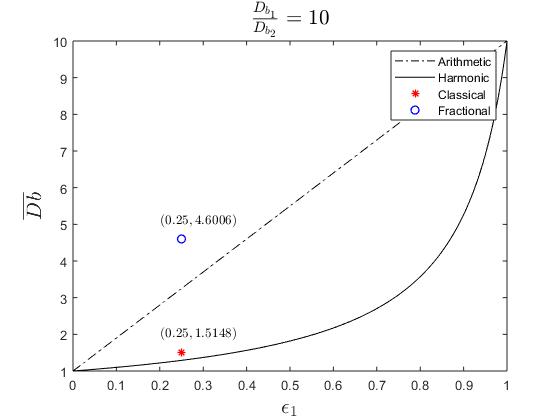}}
\end{minipage}}
\subfloat[Morphology 2: $\frac{D_{b_1}}{D_{b_2}}=0.1$]{
\label{Fig8d}
\begin{minipage}[t]{0.35\textwidth}
\centering
\scalebox{0.35}[0.35]{\includegraphics{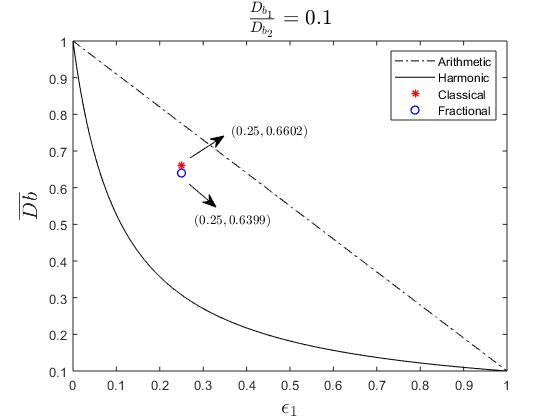}}
\end{minipage}}\\
\subfloat[Morphology 3: $\frac{D_{b_1}}{D_{b_2}}=10$]{
\label{Fig8e}
\begin{minipage}[t]{0.35\textwidth}
\centering
\scalebox{0.35}[0.35]{\includegraphics{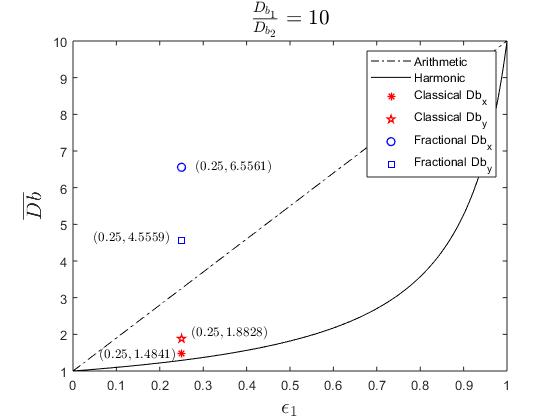}}
\end{minipage}}
\subfloat[Morphology 3: $\frac{D_{b_1}}{D_{b_2}}=0.1$]{
\label{Fig8f}
\begin{minipage}[t]{0.35\textwidth}
\centering
\scalebox{0.35}[0.35]{\includegraphics{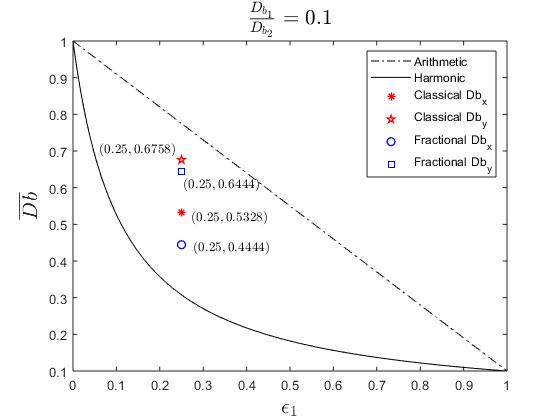}}
\end{minipage}}
\caption{The equivalent diffusivity for three morphologies with different diffusivity ratios $\frac{D_{b_1}}{D_{b_2}}$, in which we fix $D_{b_2}=1$.}
\label{Fig8}
\end{figure}

Figures \ref{Fig9a} and \ref{Fig9b} display the evolution of the equivalent diffusivity $\overline{D_{b_x}}$ for different $\alpha_1$ with increasing time and $\alpha_2=0$ under different diffusivity ratios $\frac{D_{b_1}}{D_{b_2}}$. An interesting observation from these figures is that for this two-layered composite material, the equivalent diffusivity tends to the diffusivity of the material with memory (i.e. $\overline{D_{b_x}}\rightarrow D_{b_1}=10$ when $D_{b_1}=10$ and $\overline{D_{b_x}}\rightarrow D_{b_1}=0.1$ when $D_{b_1}=0.1$), which is different to the classical case ($\overline{D_{b_x}}=1.2903$ or $0.3077$) (see Figures \ref{Fig8a} and \ref{Fig8b}). Another finding is that the fractional-order $\alpha_1$ has a significant impact on the iterative process and it delays the diffusion to reach its steady-state (at least $t>10^4s$) significantly compared to the time for the classical case ($t<1s$). We also notice that the final equivalent diffusivity $\overline{D_{b_x}}$ tends to the same constant, which is independent of the fractional-order $\alpha_1$ for a fixed $\alpha_2=0$.
\begin{table}[h]
\begin{center}
\caption{The relative error of the numerical results with the theoretical equivalent diffusivity in \cite{Szymk} for Morphology 1.}
\label{tab3}
\begin{tabular}{cccccccc}
\toprule
$D_{b_1}$  & $D_{b_2}$      & $\overline{D_{b_x}}$ & Theoretical & Relative error & $\overline{D_{b_y}}$ & Theoretical & Relative error\\
\midrule
10    & 1 &  1.2903  & 1.290 &  2.50E-04  &3.2500& 3.250 & 0\\
0.1   & 1 &  0.3077  & 0.309 &  4.23E-03 & 0.7750 & 0.775 & 0 \\
\bottomrule
\end{tabular}
\end{center}
\end{table}
\vspace{-6mm}
\vspace{-3mm}
We believe these special properties for the composite material with memory are due to the impact of the fractional operators. In addition, it appears that the fractional-order $\alpha_1$ has no influence on the equivalent diffusivity $\overline{D_{b_y}}$, which may be due to the homogeneity in the $y$ direction (see Figures \ref{Fig8a} and \ref{Fig8b}). Now, we fix $D_{b_2}=1$ and decrease the value $D_{b_1}$ to observe the evolution of the equivalent diffusivity for the classical and fractional cases ($\alpha_1=0.9$) in Figure \ref{Fig11a} using log-log axis scales. It can be seen that the equivalent diffusivity decreases with decreasing $D_{b_1}$ for both cases but they never meet (see Figure \ref{Fig11a}). Figures \ref{Fig12a} and \ref{Fig12b} depict the steady-state homogenised solutions for the classical two-layered problem and the time-fractional two-layered problem. An interesting outcome is that a very different diffusion profile is obtained for the case with memory. Again it can be concluded that the fractional-order has a significant effect on the diffusion process.\medskip

\noindent {\bf{Morphology 2}}: \emph{Circular inclusion}.\medskip

\noindent At first, we consider the case $\alpha_1=\alpha_2=0$. Figures \ref{Fig8c} and \ref{Fig8d} illustrate the equivalent diffusivity of the binary medium under different diffusivity ratios $\frac{D_{b_1}}{D_{b_2}}$. We determined that when $\frac{D_{b_1}}{D_{b_2}}=10$, $\overline{D_b}=1.5148$ and when $\frac{D_{b_1}}{D_{b_2}}=0.1$, $\overline{D_b}=0.6602$.
\vspace{-3mm}
The equivalent diffusivity for this model is between the harmonic average and arithmetic average, which is different to the two-layered problem considered on binary medium 1. In Table \ref{tab4}, the numerical equivalent diffusivity is compared with the theoretical value given in \cite{Szymk} using the initial condition $u(x,y,0)=v(x,y,0)=1$. We can see the numerical equivalent diffusivity agrees very well with the theoretical value.

Next, we observe the case when $\Omega_1$ has memory, i.e., $\alpha_1\neq 0$, $\alpha_2=0$ and the evolution of the equivalent diffusivity $\overline{D_b}$ for different $\alpha_1$ with increasing time under different diffusivity ratios $\frac{D_{b_1}}{D_{b_2}}$ is depicted in Figures \ref{Fig9c} and \ref{Fig9d}, respectively. When $\frac{D_{b_1}}{D_{b_2}}=10$, $\overline{D_b}\approx4.6006$ and when $\frac{D_{b_1}}{D_{b_2}}=0.1$, $\overline{D_b}\approx0.6399$, in which the $\overline{D_b}$ value is not bounded between the harmonic average and arithmetic average when $\frac{D_{b_1}}{D_{b_2}}=10$ (see Figures \ref{Fig8c} and \ref{Fig8d}). Also, the equivalent diffusivity $\overline{D_b}$ is independent of $\alpha_1$ with fixed $\alpha_2=0$.
\begin{table}[h]
\begin{center}
\caption{The relative error of the numerical results with the theoretical equivalent diffusivity in \cite{Szymk} for Morphology 2.}
\label{tab4}
\begin{tabular}{ccccc}
\toprule
$D_{b_1}$  & $D_{b_2}$      & $\overline{D_b}$ & Theoretical & Relative error \\
\midrule
10    & 1 &  1.5148  & 1.5200 &  3.43E-03  \\
0.1   & 1 &  0.6602  & 0.6590 &  1.85E-03  \\
\bottomrule
\end{tabular}
\end{center}
\end{table}
\vspace{-5mm}
\begin{figure}[H]
\begin{center}
\scalebox{0.3}[0.3]{\includegraphics{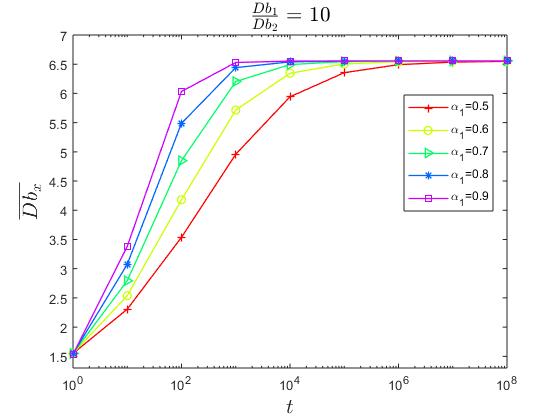}}~
\scalebox{0.3}[0.3]{\includegraphics{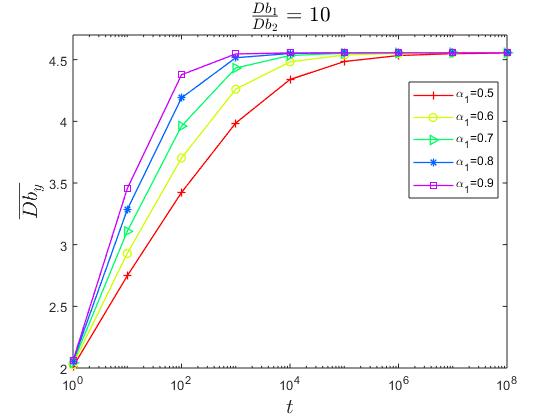}}~
\scalebox{0.3}[0.3]{\includegraphics{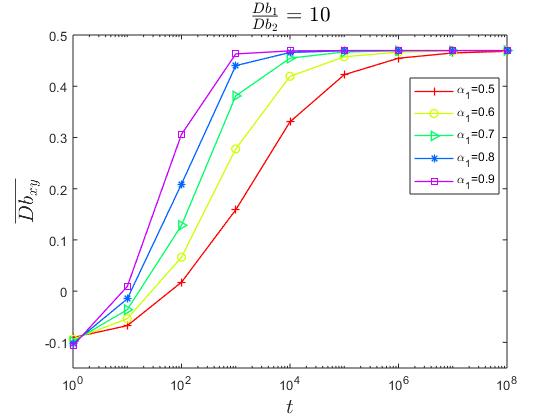}}\\
\scalebox{0.3}[0.3]{\includegraphics{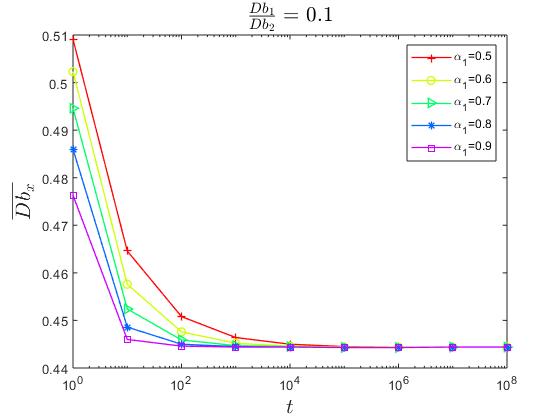}}~
\scalebox{0.3}[0.3]{\includegraphics{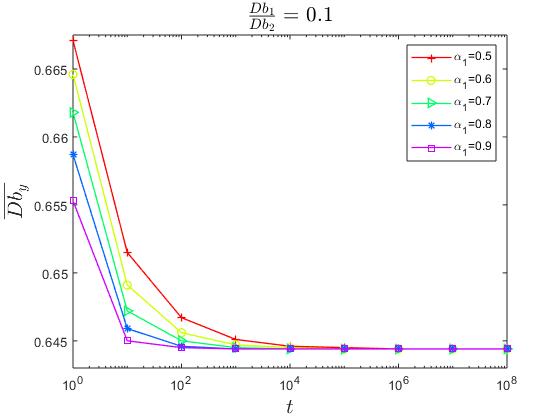}}~
\scalebox{0.3}[0.3]{\includegraphics{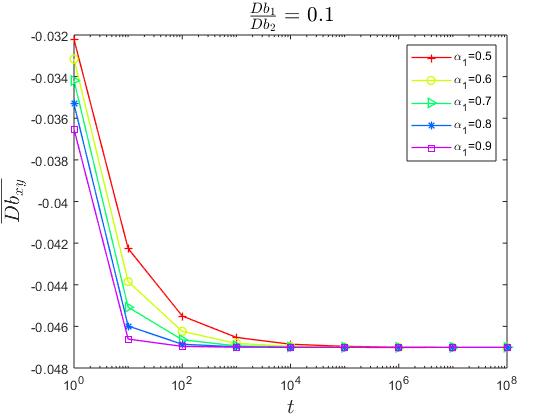}}
 \caption{The evolution of the equivalent diffusivity $\overline{D_b}$ for different $\alpha_1$ for Morphology 3 with increasing time under different medium diffusivity ratios $\frac{D_{b_1}}{D_{b_2}}$, where we fix the parameters $D_{b_2}=1$ and $\alpha_2=0$.}
\label{Fig10}
\end{center}
\end{figure}
\vspace{-5mm}
The effects of the fraction-order $\alpha_1$ on the diffusion are similar to that reported for Morphology 1. Furthermore, we study the effects of the diffusivity ratio $\frac{D_{b_1}}{D_{b_2}}$ for fixed $D_{b_2}=1$ on the equivalent diffusivity and the related results are displayed in Figure \ref{Fig11b}. It can be seen that, with a decreasing diffusivity for the memory part $D_{b_1}$, the memory effect is weakened and the equivalent diffusivity $\overline{D_b}$ coincides with the value for the classical case when $D_{b_1}$ is small ($D_{b_1}\leq 10^{-6}$), which tends to a constant $\overline{D_b}\rightarrow 0.5999$. Finally, we give a comparison of the steady-state solution profile for Morphology 2 with and without memory in Figures \ref{Fig12c} and \ref{Fig12d}. It can be observed that the memory of material $\Omega_1$ has a significant impact on the steady-state solution of the model.
\begin{figure}[H]
\centering
\subfloat[Morphology 1]{
\label{Fig11a}
\begin{minipage}[t]{0.5\textwidth}
\centering
\scalebox{0.3}[0.3]{\includegraphics{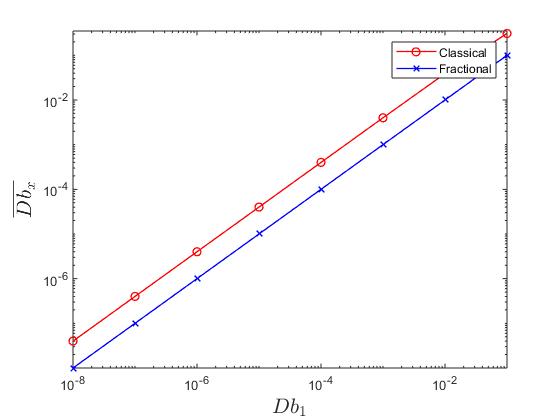}}
\end{minipage}}
\subfloat[Morphology 2]{
\label{Fig11b}
\begin{minipage}[t]{0.5\textwidth}
\centering
\scalebox{0.3}[0.3]{\includegraphics{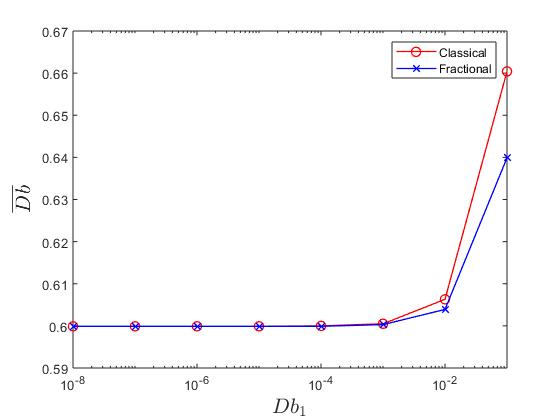}}
\end{minipage}}
\caption{The effects of the diffusivity ratio $\frac{D_{b_1}}{D_{b_2}}$ for fixed $D_{b_2}=1$ on the equivalent diffusivity, where we use the parameter $\alpha_1=0.9$ for the fractional case. The equivalent diffusivity is presented using a log-log axis (a) and semi-log axis (b), respectively.}
\label{Fig11}
\end{figure}
\vspace{-3mm}
\noindent {\bf{Morphology 3}}: \emph{L-shaped inclusion}.\medskip

\noindent Again, we first consider the classical $\alpha_1=\alpha_2=0$ case. Figures \ref{Fig8e} and \ref{Fig8f} illustrate the equivalent diffusivity for the problem on binary medium 3 under different medium diffusivity ratios $\frac{D_{b_1}}{D_{b_2}}$ with initial condition is $u(x,y,0)=v(x,y,0)=1$. The relative error between the numerical equivalent diffusivity and the theoretical diffusivity is presented in Table \ref{tab5}.
\begin{table}[h]
\begin{center}
\caption{The relative error of the numerical results with the theoretical equivalent diffusivity in \cite{Szymk} for Morphology 3.}
\label{tab5}
\begin{tabular}{ccccc}
\toprule
$D_{b_1}$  & $D_{b_2}$      & $\overline{D_{b_x}}$ & Theoretical & Relative error \\
\midrule
10    & 1 &  1.4841  & 1.4800 &  2.75E-03  \\
0.1   & 1 &  0.5328  & 0.5330 &  3.72E-04  \\
\midrule
$D_{b_1}$  & $D_{b_2}$      & $\overline{D_{b_y}}$ & Theoretical & Relative error \\
\midrule
10    & 1 &  1.8828  & 1.8800 &  1.49E-03  \\
0.1   & 1 &  0.6758  & 0.6750 &  1.15E-03  \\
\midrule
$D_{b_1}$  & $D_{b_2}$      & $\overline{D_{b_{xy}}}$ & Theoretical & Relative error \\
\midrule
10    & 1 &  -7.9601E-02  & -0.0796 &  1.80E-05  \\
0.1   & 1 &  -2.8599E-02  & -0.0286 &  4.51E-05  \\
\bottomrule
\end{tabular}
\end{center}
\end{table}
\vspace{-3mm}

When $\Omega_1$ has memory ($\alpha_1\neq 0$, $\alpha_2=0$), the evolution of the equivalent diffusivity $\overline{D_b}$ for different $\alpha_1$ with increasing time under different medium diffusivity ratios $\frac{D_{b_1}}{D_{b_2}}$ is displayed in Figure \ref{Fig10}. We see that the equivalent diffusivity tends to a fixed value at steady-state, which is consistent with the findings reported for Morphologies 1 and 2. When $\frac{D_{b_1}}{D_{b_2}}=10$, $\overline{D_{b_x}}\to 6.5561$, $\overline{D_{b_y}}\to 4.5559$, $\overline{D_{b_{xy}}}\to 0.46969$ and when $\frac{D_{b_1}}{D_{b_2}}=0.1$, $\overline{D_{b_x}}\to 0.4444$, $\overline{D_{b_y}}\to 0.6444$, $\overline{D_{b_{xy}}}\to -0.047007$. For this morphology, the equivalent diffusivity values $\overline{D_{b_x}}$ and $\overline{D_{b_y}}$ for the classical case all lie between the harmonic average and arithmetic average, while when $\frac{D_{b_1}}{D_{b_2}}=10$, the values for the fractional case are out of this bound (see Figures \ref{Fig8e} and \ref{Fig8f}). An interesting finding is that the sign of the cross-diffusion component $\overline{D_{b_{xy}}}$ has changed in the whole iterative process when $\frac{D_{b_1}}{D_{b_2}}=10$ (see Fig10). Finally, the comparison of the steady-state in the $x$-direction and $y$-direction between the classical case and fractional case is shown in Figures \ref{Fig12d} to \ref{Fig12h}, from which some clear differences can be observed in the steady-state solution.

In summary, the equivalent diffusivity for the classical two-layered medium can be bounded by the harmonic average and arithmetic average. However, the equivalent diffusivity for the binary medium with memory effect is not bounded by these limits. For the fractional homogenisation problem, the equivalent diffusivity value tends to a fixed value after a long time iteration. In addition, the equivalent diffusivity deceases and recovers the classical case as the diffusivity coefficient of the memory part is decreased to zero.
\begin{figure}[H]
\centering
\subfloat[Morphology 1]{
\label{Fig12a}
\begin{minipage}[t]{0.25\textwidth}
\centering
\scalebox{0.22}[0.22]{\includegraphics{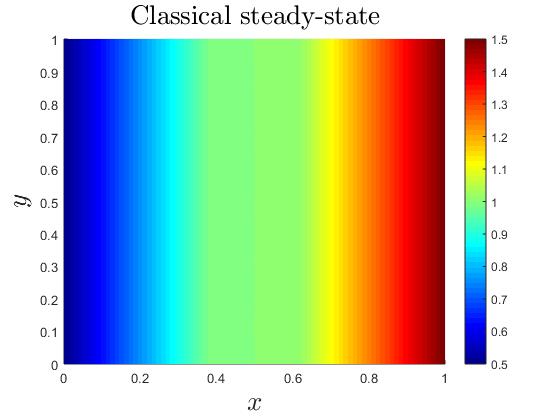}}
\end{minipage}}
\subfloat[Morphology 1]{
\label{Fig12b}
\begin{minipage}[t]{0.25\textwidth}
\centering
\scalebox{0.22}[0.22]{\includegraphics{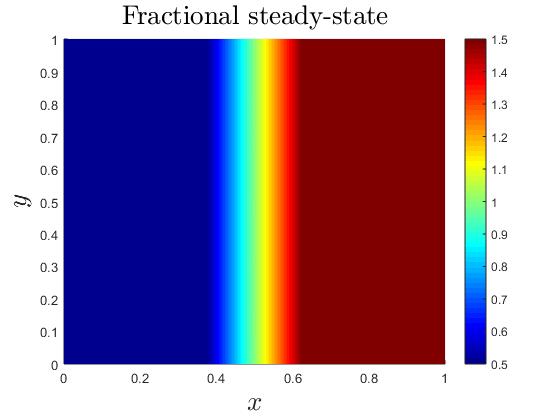}}
\end{minipage}}
\subfloat[Morphology 2]{
\label{Fig12c}
\begin{minipage}[t]{0.25\textwidth}
\centering
\scalebox{0.22}[0.22]{\includegraphics{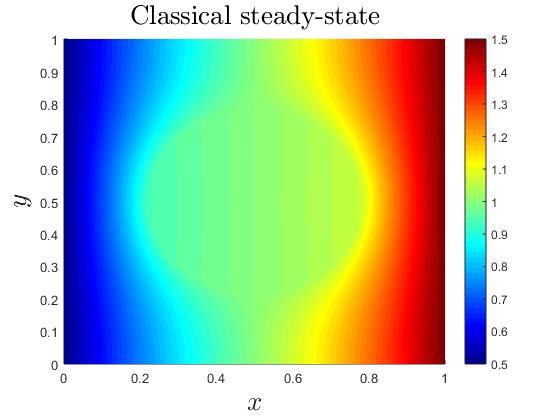}}
\end{minipage}}
\subfloat[Morphology 2]{
\label{Fig12d}
\begin{minipage}[t]{0.25\textwidth}
\centering
\scalebox{0.22}[0.22]{\includegraphics{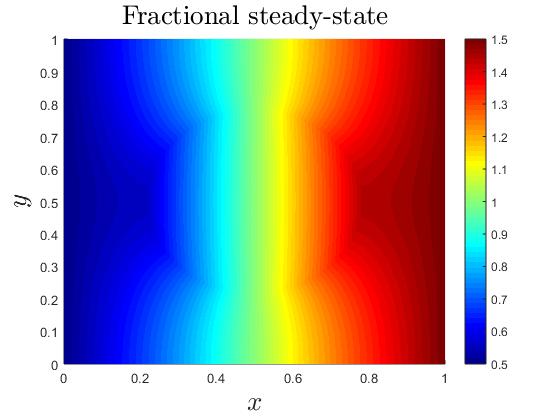}}
\end{minipage}}\\
\subfloat[Morphology 3]{
\label{Fig12e}
\begin{minipage}[t]{0.25\textwidth}
\centering
\scalebox{0.22}[0.22]{\includegraphics{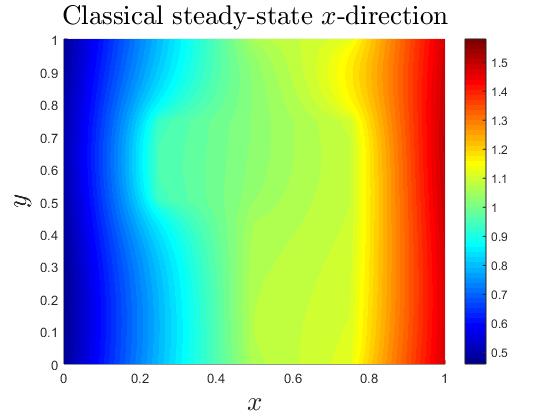}}
\end{minipage}}
\subfloat[Morphology 3]{
\label{Fig12f}
\begin{minipage}[t]{0.25\textwidth}
\centering
\scalebox{0.22}[0.22]{\includegraphics{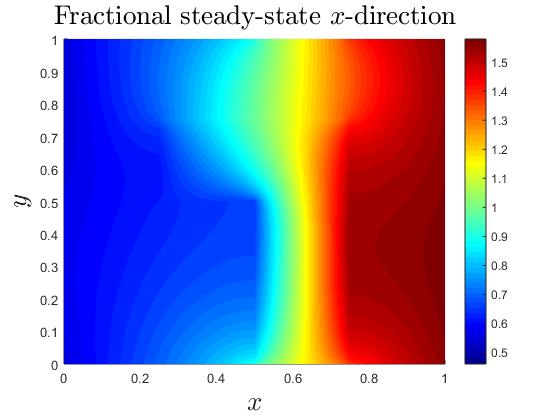}}
\end{minipage}}
\subfloat[Morphology 3]{
\label{Fig12g}
\begin{minipage}[t]{0.25\textwidth}
\centering
\scalebox{0.22}[0.22]{\includegraphics{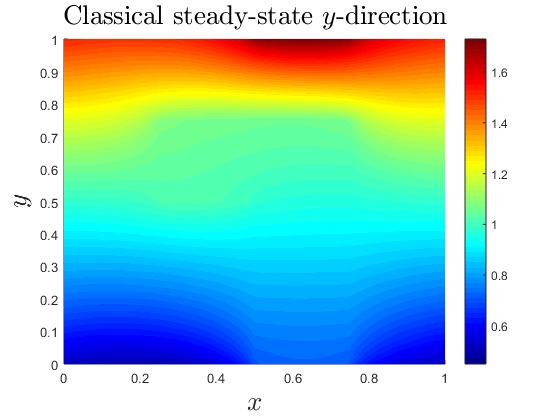}}
\end{minipage}}
\subfloat[Morphology 3]{
\label{Fig12h}
\begin{minipage}[t]{0.25\textwidth}
\centering
\scalebox{0.22}[0.22]{\includegraphics{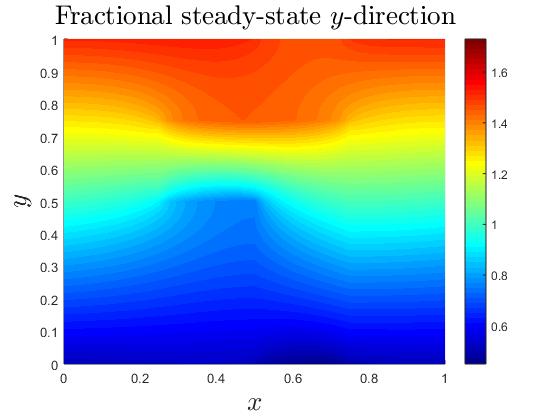}}
\end{minipage}}
\caption{The steady-state for classical and time-fractional case of the three different kinds of two-layered problems with $\frac{D_{b_1}}{D_{b_2}}=10$. }
\label{Fig12}
\end{figure}
\vspace{-5mm}
\section{Application of the model to estimating the bound water diffusivity of wood}\label{sec6}

\begin{figure}[H]
\centering
\subfloat[Spruce earlywood]{
\label{Fig13a}
\begin{minipage}[t]{0.8\textwidth}
\centering
{\includegraphics[width=0.6\textwidth]{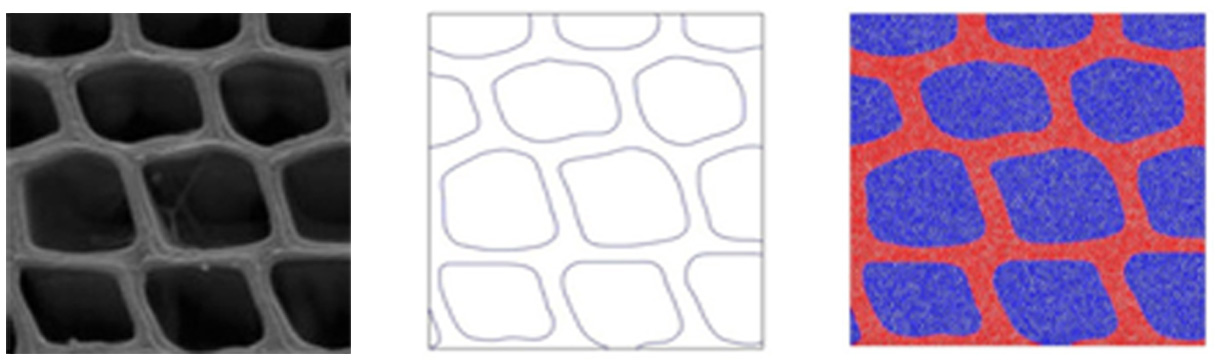}}
\end{minipage}
}\\
\subfloat[Spruce latewood]{
\label{Fig13b}
\begin{minipage}[t]{0.8\textwidth}
\centering
{\includegraphics[width=0.6\textwidth]{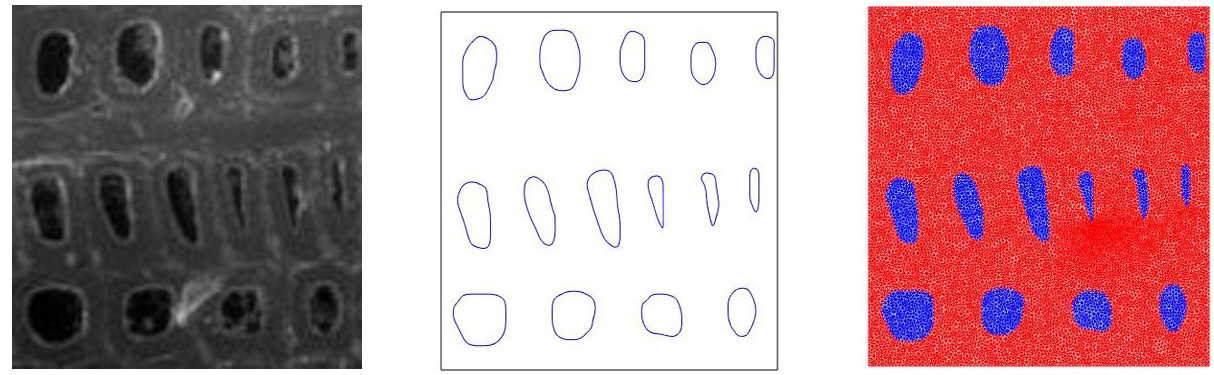}}
\end{minipage}
}\\
\subfloat[Australian hardwood]{
\label{Fig13c}
\begin{minipage}[t]{0.8\textwidth}
\centering
{\includegraphics[width=0.6\textwidth]{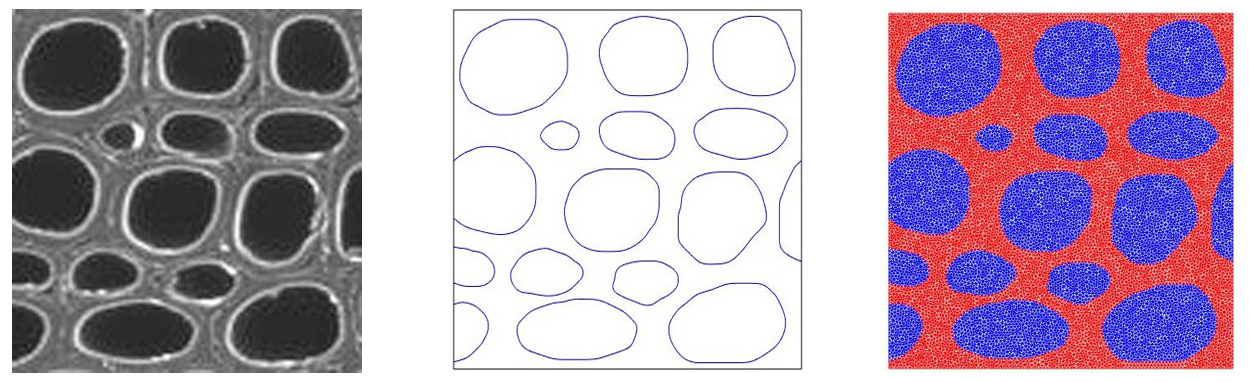}}
\end{minipage}
}
\caption{The anatomical images of two different species of wood cells ($100\mu m\times100\mu m$) and their boundary contour and triangulation, in which the volumetric fraction of solid part ($\Omega_2$) is approximately (a) 0.3427, (b) 0.8445, (c) 0.4566, respectively.}
\label{Fig13}
\end{figure}

We now use the time-fractional transport model to perform computational homogenisation using different microscale images of Spruce earlywood, Spruce latewood and an Australian hardwood (blackbutt) to estimate the bound water diffusivity when the cellular structure exhibits anomalous diffusion or memory effects. This type of behaviour has been observed experimentally in \cite{Turner11,Olek,Olek2,Perr}. The ESEM imaging of these wood species were performed at the Laboratoire G\'{e}nie des Proc\'{e}d\'{e}s et Mat\'{e}riaux, Centrale--Sup\'{e}lec in France and are shown in Figure \ref{Fig13}. Recent experimental work highlighted that traditional diffusion equations cannot adequately describe the absorption of water in a cell wall \cite{Olek,Perr}. These effects are exacerbated even further for a modified cell wall occurring due to heat treatment \cite{Olek2}. The reduction in diffusivity is thought to be due to the degradation of hemicelluloses. The thermal treatment leads to a change in the mechanism of attaining hygroscopic equilibrium by the modified wood. For these applications, traditional conservation laws prove inadequate for describing the underlying physical processes, while models based on fractional operators succeed \cite{Turner11}.

Fractional-order derivatives provide excellent alternatives to their classical counterparts for such applications by interpolating between the integer orders of differential equations to capture nonlocal relations in time using power-law memory kernels. A major attraction is that the dissipative and dispersive properties observed in experimental data are representative of a variety of different forms of generalised constitutive laws. Diffusion involving molecular relaxation in the complex porous microstructure induces a second time constant in the macroscopic physical phenomena that gives rise to `sub-diffusive' transport behaviour (see chapter 17 \cite{Klages}), which motivates why we have used a fractional-in-time derivative for the generalised conservation law. For the computational homogenisations performed over the given unit cell pore structures, the following modified transient diffusion equation can be used to model the transport of water within the wood cell, which is an extension of the wood drying problem discussed in \cite{Carr13,Foy17}
\begin{align*}
\frac{\partial \Psi_w}{\partial t}+{_0^RD^{\alpha}_{t}}(\nabla \cdot {\bf{Q}}_w)=0.
\end{align*}
Introducing the indicator variable
\begin{equation*}
\chi=\left\{\begin{array}{ll}
1,& (x,y)\in C^{(g)},\\
0,& (x,y)\in C^{(s)},
\end{array}\right.
\end{equation*}
then we can define the fractional indexes and the mass fluxes as $\alpha=(1-\chi)\alpha_1+\chi\alpha_2$, ${\bf{Q}}_w=(1-\chi){\bf{Q}}_w^{(s)}+\chi {\bf{Q}}_w^{(g)}$, where
${\bf{Q}}_w^{(s)}=-\rho_sD_b\nabla X_m$,${\bf{Q}}_w^{(g)}=-\frac{\rho_gD_v}{1-\omega_v}\nabla\omega_v$, where $X_m$ is the cell wall bound water moisture content and $\rho_v$ is the vapour density. The conserved quantity is defined as $\Psi_w=(1-\chi)\rho_sX_m+\chi\rho_v$. The quasi-periodic boundary conditions are
\begin{align*}
&X_m(L_x,y)=X_m(0,y)+\frac{1}{L_x}\frac{\partial X_B}{\partial x},~ 0<y<L_y,\\
&X_m(x,L_y)=X_m(x,0)+\frac{1}{L_y}\frac{\partial X_B}{\partial y},~ 0<x<L_x,\\
&{_0^RD^{\alpha}_{t}}[{\bf{Q}}_w(0,y)\cdot {\bf{n}}_{\partial C}]={_0^RD^{\alpha}_{t}}[{\bf{Q}}_w(L_x,y)\cdot {\bf{n}}_{\partial C}],~ 0<y<L_y,\\
&{_0^RD^{\alpha}_{t}}[{\bf{Q}}_w(x,0)\cdot {\bf{n}}_{\partial C}]={_0^RD^{\alpha}_{t}}[{\bf{Q}}_w(x,L_y)\cdot {\bf{n}}_{\partial C}],~ 0<x<L_x,
\end{align*}
where ${\bf{n}}_{\partial C}$ is the unit vector normal to $\partial C$ outward to $C$, $\frac{\partial X_B}{\partial x}$ and $\frac{\partial X_B}{\partial y}$ are the gradients imposed over the cell, respectively. The closure conditions  for the problem are
$P_a+P_v=P_{atm}$, $\rho_v=\frac{P_vM_v}{R(T+273.15)}$, $\rho_a=\frac{P_aM_a}{R(T+273.15)}$, where the parameters are calculated using the same formulae in \cite{Carr13} at $T=20$\textcelsius.
\vspace{-3mm}

Since wood has a periodic structure \cite{Perre16}, a representative elementary volume (unit cell), needs to be chosen. We used ESEM images of the real wood cellular structure by extracting the boundary of the pores and then a digital representation of the cell was formed. Furthermore, we use the mesh generator \emph{Gmsh} \cite{Geuzaine} to perform the triangulation, which can divide the cell domain $[0,L_x]\times[0,L_y]$ into two sub-sets: the cell lumens ($\Omega_1$ or $C^{(g)}$) and the solid phase ($\Omega_2$ or $C^{(s)}$). Here, we refer the $x$-direction as the radial direction and refer the $y$-direction as the tangential direction. We denote the harmonic average as the series value and the arithmetic average as the parallel value to be consistent with \cite{Perre16}. For the treatment of the periodic boundary conditions, one can refer to \cite{Carr13}. Here we will consider the problem on two different species of wood cells: softwood Spruce earlywood and latewood and Australian hardwood (see Figure \ref{Fig13}).
\begin{figure}[H]
\centering
\subfloat[Spruce earlywood]{
\label{Fig14a}
\begin{minipage}[t]{0.33\textwidth}
\centering
\scalebox{0.32}[0.32]{\includegraphics{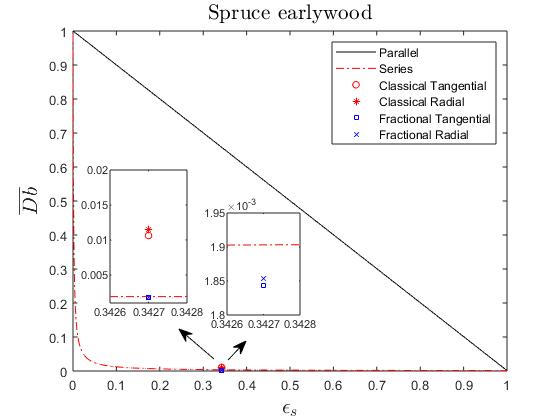}}
\end{minipage}
}
\subfloat[Spruce latewood]{
\label{Fig14b}
\begin{minipage}[t]{0.33\textwidth}
\centering
\scalebox{0.32}[0.32]{\includegraphics{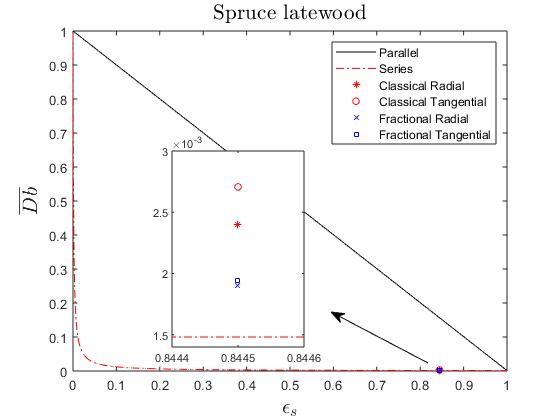}}
\end{minipage}
}
\subfloat[Australian hardwood]{
\label{Fig14c}
\begin{minipage}[t]{0.33\textwidth}
\centering
\scalebox{0.32}[0.32]{\includegraphics{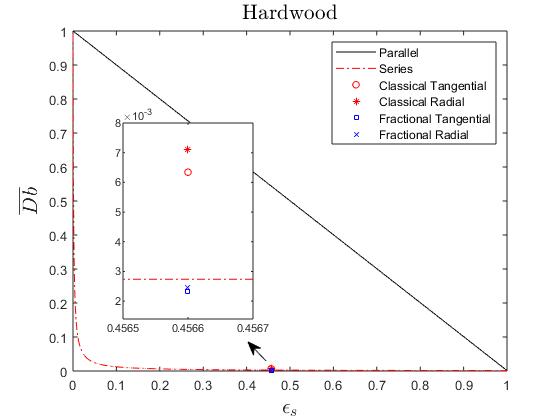}}
\end{minipage}
}
\caption{The equivalent diffusivity for three different species of wood cells. Note that all the values have been scaled through division by $D_v$.}
\label{Fig14}
\end{figure}
\vspace{-3mm}

Figure \ref{Fig14a} shows the dimensionless equivalent diffusivity of Spruce earlywood for two different cases: the solid parts without memory (classical case) and with memory (fractional case). Since the solid cell walls have a low diffusivity of water, the equivalent diffusivity value is very close to the series model. As the cell walls in earlywood are aligned with the radial direction, the value in the radial direction is larger than in the tangential direction. For the classical case, the radial value ($\overline{D_{b_x}}\approx 1.0628\times 10^{-2}$) and tangential value ($\overline{D_{b_y}}\approx 1.1517\times 10^{-2}$) are bounded by the series and parallel values, which agrees with the results presented in \cite{Foy17,Perre16}.  While for the fractional case (where the solid part has memory), the radial value ($\overline{D_{b_x}}\approx 1.8437\times 10^{-3}$) and tangential value ($\overline{D_{b_y}}\approx 1.8534\times 10^{-3}$) are not bounded by the series and parallel values, being smaller than the classical case.

In the following discussion, we consider the equivalent diffusivity of Spruce latewood cell problem. Figure \ref{Fig14b} displays the dimensionless equivalent diffusivity of Spruce latewood for two different cases: the solid parts without memory (classical case) and with memory (fractional case). In contrast to Spruce earlywood, the tracheids in latewood are more flattened in the tangential direction, which blocks the mass flux in the radial direction. An inverted anisotropic ratio is observed. For the classical case, the radial value is $\overline{D_{b_x}}\approx 2.7068\times 10^{-3}$ and tangential value is $\overline{D_{b_y}}\approx 2.3991\times 10^{-3}$. For the fractional case (the solid part has memory), the radial value is $\overline{D_{b_x}}\approx 1.9415\times 10^{-3}$ and tangential value is $\overline{D_{b_y}}\approx 1.9009\times 10^{-3}$, both of which are smaller than the classical case. Next, we consider the Australian hardwood case. The dimensionless equivalent diffusivity of hardwood for the classical case and the fractional case is presented in Figure \ref{Fig14c}, from which a similar result can be seen.  For the classical case, the radial value ($\overline{D_{b_x}}\approx 7.0923\times 10^{-3}$) and tangential value ($\overline{D_{b_y}}\approx 6.3437\times 10^{-3}$) are bounded by the series and parallel values. While for the fractional case (the solid part has memory), the radial value ($\overline{D_{b_x}}\approx 2.4440\times 10^{-3}$) and tangential value ($\overline{D_{b_y}}\approx 2.3370\times 10^{-3}$) are not bounded by the series and parallel values.
\begin{figure}[H]
\centering
\subfloat[Spruce earlywood]{
\label{Fig15a}
\begin{minipage}[t]{0.3\textwidth}
\centering
\scalebox{0.4}[0.4]{\includegraphics{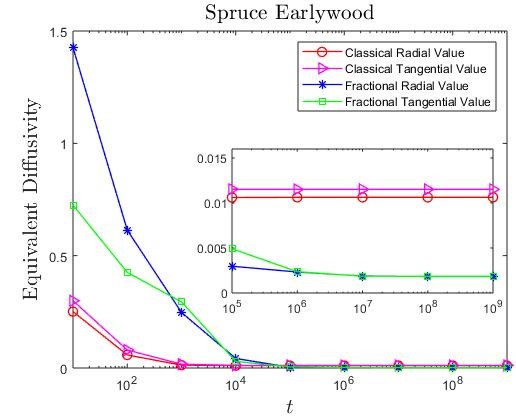}}
\end{minipage}
}
\subfloat[Spruce latewood]{
\label{Fig15b}
\begin{minipage}[t]{0.3\textwidth}
\centering
\scalebox{0.4}[0.4]{\includegraphics{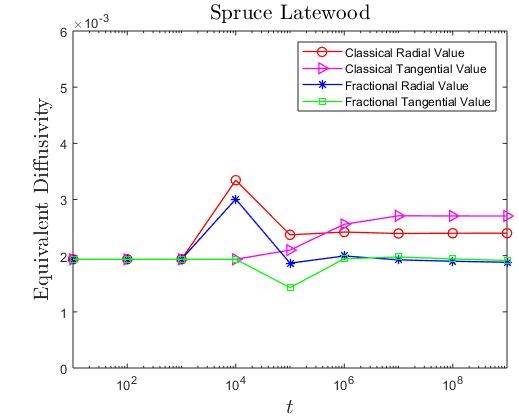}}
\end{minipage}
}
\subfloat[Australian hardwood]{
\label{Fig15c}
\begin{minipage}[t]{0.3\textwidth}
\centering
\scalebox{0.4}[0.4]{\includegraphics{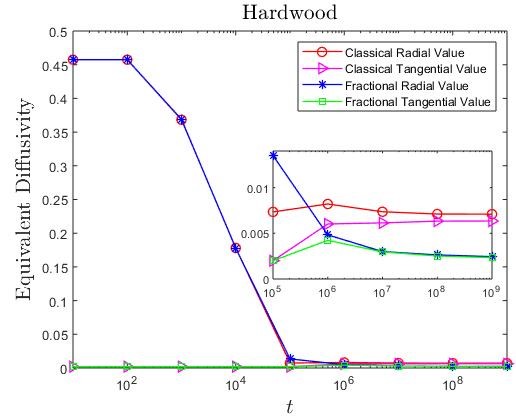}}
\end{minipage}
}
\caption{The evolution of the equivalent diffusivity for three different species of wood cells. Note that all the values have been scaled through division by $D_v$.}
\label{Fig15}
\end{figure}

Finally, Figure \ref{Fig15} compares the evolution of the equivalent diffusivity for the three different species of wood cells, namely Spruce earlywood, Spruce latewood and Australian hardwood. It is can be seen that, in accordance with the findings reporting in Section \ref{sec5}, the effective bound water diffusivities slowly approach the values of the solid phase (that undergoes molecular relaxation) diffusivities as the homogenised cell problem reaches steady state. Again, this process can take considerable time for the fractional model and tends to be much more rapid for classical diffusion.

In summary, for Spruce earlywood without memory, the tangential value is larger than the radial value while an inverted anisotropic ratio is observed for Spruce latewood. For Australian hardwood, the tangential and radial values are between those of Spruce earlywood and Spruce latewood, which is due to the cell structure. For both species of wood, the tangential and radial values are bounded by the series and parallel values.  However, the tangential and radial values are not bounded by these values when the wood has memory. In addition, the tangential and radial values for wood with memory are smaller than those for wood without memory, which we believe is due to the memory effect. To conclude, for the homogenisation of the time-fractional transport problem, the fractional-order indices have a significant impact on the mass transfer and the equivalent diffusivity appears to be dominated by the lignocellulosic material with memory effects.

\section{Conclusions}\label{sec7}

In this paper, we considered a two-dimensional time-fractional subdiffusion equation, in which the unstructured mesh control volume method is applied. Our chosen numerical examples showed that the method was stable and effective. As an application, we successfully simulated a time-fractional transport model in a binary medium consisting of regular and irregular inclusions and derived a semi-analytical solution for a class of two-layered problems subjected to quasi-periodic boundary conditions. An important contribution was the extension of the classical homogenisation theory to accommodate the new framework and to show that the effective diffusivity tensor can be computed once the cell problems reach steady state. We found for all morphologies considered that the effective diffusivity slowly converged to the diffusivity of the material in the unit cell with memory, which was a very different finding to that observed for homogenised diffusivity parameters for standard materials.

In order to harness the main findings of our work, we postulate that the time evolutionary behaviour of the effective parameters must be used in the homogenised macroscopic model until steady-state in the unit cell is achieved, at which time the homogenised effective parameter is dominated by the diffusivity of the material exhibiting memory effects. Furthermore, the fractional index used in the macroscopic time-fractional equation must be averaged over the unit cell. In our future research, we plan to investigate this homogenisation strategy in more detail and verify it using experimental data. We will also use the time-fractional multi-scale problem to simulate more complicated heterogeneous systems.

\section*{Acknowledgment}
This work was supported financially by a visiting Professorial Fellowship that enabled Turner to work at the Universit\'{e} Paris-Saclay, France for a period of three months in 2017-18. We acknowledge the financial support for this research received through the Australian Research Council Discovery Grant DP150103675.

\end{document}